\documentclass[journal]{IEEEtran}
\usepackage{amsmath,amsthm,amssymb}
\usepackage{cite}
\usepackage[hidelinks]{hyperref} 
\usepackage[caption=false,font=footnotesize]{subfig}
\usepackage{psfrag}
\usepackage{stfloats}
\fnbelowfloat
\usepackage{array}
\usepackage{algorithmicx}
\usepackage{algorithm}
\newtheorem{mylemma}{Lemma}
\newtheorem{myproposition}{Proposition}
\newtheorem{mytheorem}{Theorem}

\usepackage{ifpdf}
 \ifpdf
 \else
 \fi

\ifCLASSINFOpdf
 \usepackage[pdftex]{graphicx}
 \graphicspath{{../pdf/}{../jpeg/}}
 \DeclareGraphicsExtensions{.pdf,.jpeg,.png}
\else
 \usepackage[dvips]{graphicx}
 \graphicspath{ {Figures/} }
 \DeclareGraphicsExtensions{.eps}
\fi

\newtheoremstyle{mytheoremstyle}{0pt}{0pt}{\itshape}{}{\bfseries}{.}{.5em}{} 
\theoremstyle{mytheoremstyle}

\makeatletter
\def\thmhead@plain#1#2#3{%
	\thmname{#1}\thmnumber{\@ifnotempty{#1}{ }\@upn{#2}}%
	\thmnote{ {\the\thm@notefont#3}}}
\let\thmhead\thmhead@plain
\def\thm@space@setup{\thm@preskip=0pt
	\thm@postskip=0pt}
\makeatother

\usepackage{xpatch}
\makeatletter
\xpatchcmd{\proof}{\@addpunct{.}}{\@addpunct{:}}{}{} 
\xpatchcmd{\proof}{\hskip\labelsep}{\hskip5\labelsep}{}{} 
\xpatchcmd{\proof}{\topsep6\p@\@plus6\p@\relax}{}{}{}
\makeatother

\usepackage[noend]{algpseudocode}
\begin{document}
\title{Design of PAR-Constrained Sequences for MIMO Channel Estimation via Majorization-Minimization}
\author{Zhongju~Wang,
        Prabhu~Babu,
        and~Daniel~P.~Palomar,~\IEEEmembership{Fellow,~IEEE}
\thanks{This work was supported by the Hong Kong RGC 16206315 research
	grant. Zhongju Wang and Daniel P. Palomar are with the Hong Kong
	University of Science and Technology (HKUST), Hong Kong. E-mail: {\{zwangaq, palomar\}}@ust.hk. Prabhu Babu was with the Hong Kong University of Science and Technology (HKUST), Hong Kong. He is now with CARE, IIT Delhi, Delhi, India. Email: prabhubabu@care.iitd.ac.in.}
}


\maketitle

\begin{abstract}
 PAR-constrained sequences are widely used in communication systems and radars due to various practical needs; specifically, sequences are required to be unimodular or of low peak-to-average power ratio (PAR). For unimodular sequence design, plenty of efforts have been devoted to obtaining good correlation properties. Regarding channel estimation, however, sequences of such properties do not necessarily help produce optimal estimates. Tailored unimodular sequences for the specific criterion concerned are desirable especially when the prior knowledge of the channel is taken into account as well. In this paper, we formulate the problem of optimal unimodular sequence design for minimum mean square error estimation of the channel impulse response and conditional mutual information maximization, respectively. Efficient algorithms based on the majorization-minimization framework are proposed for both problems with guaranteed convergence. As the unimodular constraint is a special case of the low PAR constraint, optimal sequences of low PAR are also considered. Numerical examples are provided to show the performance of the proposed training sequences, with the efficiency of the derived algorithms demonstrated.
\end{abstract}

\begin{IEEEkeywords}
Unimodular sequence, peak-to-average power ratio (PAR), channel estimation, majorization-minimization, minimum mean square error, conditional mutual information.
\end{IEEEkeywords}
%
\IEEEpeerreviewmaketitle

\section{Introduction}
\label{sec:intro}
\IEEEPARstart{P}{AR-constrained} sequences, such as unimodular or low peak-to-average power ratio (PAR), have many applications in both single-input single-output (SISO) and multi-input multi-output (MIMO) communication systems. For example, the $M$-ary phase-shift keying techniques allow only symbols of constant-modulus, i.e., unimodular, to be transmitted \cite{benedetto1999principles}. In MIMO radars and code-division multiple-access (CDMA) applications, the practical implementation demands from hardware, such as radio frequency power amplifiers and analog-to-digital converters, require the sequences transmitted to be unimodular or low PAR \cite{tropp2005designing,he2009designing,he2011wideband}. In this paper, we consider the design of optimal unimodular or low PAR sequences for channel estimation.

There is an extensive literature on designing single unimodular sequences with good correlation properties such that the autocorrelation of the sequence is zero at each nonzero lag. As such properties are usually difficult to achieve, metrics of ``goodness'' have been proposed instead where autocorrelation sidelobes are suppressed rather than literally set to zero, and optimization problems are thus formulated and solved with numerical algorithms \cite{stoica2009new,song2014optimal}. Specifically, the work \cite{stoica2009new} provides several cyclic algorithms (CA) for either minimizing integrated sidelobe level (ISL) or maximizing ISL-related merit factor (MF). In \cite{song2014optimal}, a computationally efficient algorithm called MISL for minimizing ISL is proposed, and it is demonstrated that MISL results in lower autocorrelation sidelobes with less computational complexity.

The good correlation property of a single unimodular sequence is also extended to MIMO systems, where multiple sequences are transmitted. The good autocorrelation is defined for each sequence as that for a single sequence. Meanwhile, good cross-correlation demands that any sequence be nearly uncorrelated with time-shifted versions of the other sequences. In \cite{he2009designing}, algorithms CA-direct (CAD) and CA-new (CAN) are developed to obtain sequence sets of low auto- and cross-correlation sidelobes. Also \cite{song2015seqencce} proposes some efficient algorithms to minimize the same metric in \cite{he2009designing}. 

The aforementioned ISL and ISL-related metrics are both alternative ways to describe the impulse-like correlation characteristics. Sequences with such properties enable matched filters at the receiver side to easily extract the signals backscattered from the range bin of interest and attenuate signals backscattered from other range bins \cite{he2009designing}. Nevertheless, matched filters take no advantage of any prior information on the channel when the unimodular low-ISL sequences are used for estimation.

The unimodular constraint is actually a special case of the low PAR constraint, which imposes how the largest amplitude of the sequence compares with its average power. The low PAR constraint, as a structural requirement, has been well studied in the design of tight frames \cite{tropp2005designing}. Although the individual vector norms of a frame could be adjusted to maximize the sum-capacity of DS-CDMA links, the optimality in terms of any performance measures was not directly considered therein. Furthermore, the algorithm they proposed is based on alternating projection that often suffers a slow convergence.

As far as channel estimation is concerned, many studies have been conducted for both frequency-flat and frequency-selective fading channels under minimum mean square error (MMSE) estimation
and conditional mutual information (CMI) maximization. Most of those obtained optimization problems, however, only address the power constraint without addressing the unimodular or low PAR constraints. In \cite{Kotecha2004,liu2007training,katselis2008training,biguesh2009optimal,bjornson2010framework}, training sequence design for flat MIMO channels is studied assuming some special structures, e.g., Kronecker product, on the prior covariance matrices of channel and noise. It is shown that the optimization problem can be reformulated as power allocations using the majorization theory, and the waterfilling solutions are obtained. Meanwhile, problems of similar formulations have also been studied in joint linear transmitter-receiver design \cite{palomar2003joint,palomar2006mimo}. To deal with arbitrarily correlated MIMO channels, some numerical algorithms based on block coordinate descent are proposed in \cite{katselis2013training,bengtsson2014algorithmic,shi2014training}. 

More related to our work is training sequence design for frequency-selective fading channels. Under a total power constraint, channel capacity is investigated for SISO channels \cite{vikalo2004capacity} and MIMO channels \cite{ma2005optimal}. Independent and identically distributed channel coefficients and noise are assumed to facilitate mathematical analysis. As a result, impulse-like sequences for both types of channel are suggested for optimal estimation. Optimal design for the MMSE channel estimation has been studied in \cite{leus2005optimal}, where the noise is assumed to be white and the channel taps are uncorrelated; however, such assumption is hardly satisfied in practice. And there is no guarantee of finding an optimal solution for an arbitrary length of training or channel correlation. More important, their results cannot be used when the unimodular or low PAR constraint is imposed on the sequences to be designed.

We formulate the problem as the design of optimal unimodular sequences based on the MMSE and the CMI. Both problems are non-convex with the bothersome unimodular constraint. Without assuming any amenable structures, e.g., Kronecker product, on the prior channel and noise covariances, the problems are also challenging even if only the power constraint is imposed. To tackle those issues, the majorization-minimization (MM) technique is employed to develop efficient algorithms. By rewriting the objective functions in a more appropriate way, majorizing/minorizing functions can be obtained for minimization/maximization objective. As a result, the original problems are solved instead by a sequence of simple problems, each of which turns out to have a closed-form solution. Convergence of our proposed algorithms is guaranteed, and an acceleration scheme is also given to improve the convergence rate. For low PAR constraints, similar problems can be formulated, and the developed algorithms need only a few modifications to be applied.

The rest of this paper is organized as follows. In Section \ref{sec:channel_problem}, the channel model is described, on which the optimal unimodular sequence design problems are formulated. In Section \ref{sec:Algorithms_for_Optimal}, derivations of algorithms for both the MMSE minimization and the CMI maximization are presented, followed by a brief analysis of convergence properties and an acceleration scheme. The optimal design under the low PAR constraints is discussed in \ref{sec:Opt_Seq_PAR}. Numerical examples are presented in Section \ref{sec:simulations}. And conclusion is then given in Section \ref{sec:conclusion}.

\textit{Notation:} Scalars are represented by italic letters. Boldface uppercase and lowercase letters denote matrices and vectors, respectively. $\mathbb{C}$ is the set of complex numbers. The identity matrix is denoted by $\mathbf{I}$ with the size implicit in the context if undeclared. The superscripts $\left(\cdot\right)^{T}$, $\left(\cdot\right)^{H}$ and $\left(\cdot\right)^{\ast}$ denote respectively transpose, conjugate transpose and complex conjugate. With  $\mathrm{vec}\left(\mathbf{X}\right)$, the vector is formed by stacking the columns of $\mathbf{X}$. The Kronecker product is denoted by $\otimes$. $\mathrm{E}\left(\cdot\right)$ takes the expectation of random variable. $\mathrm{Tr}\left(\cdot\right)$ is the trace of a matrix. $\|\cdot\|_{F}$ is Frobenius norm of a matrix.


\section{Channel Model and Problem Formulations}
\label{sec:channel_problem}
We consider a block-fading or quasistatic multi-input multi-output (MIMO) channel. Assume the number of transmit antennas and receive antennas are $N_{t}$ and $N_{r}$, respectively, and the channel impulse response is described as a length-$(K+1)$ sequence of matrices $\mathbf{H}_{0},\ldots,\mathbf{H}_{K}\in\mathbb{C}^{N_{r}\times N_{t}}$. In the training period, a length-$N$ sequence is sent through the channel from each transmit antenna or, equivalently, a length-$N_t$ vector $\mathbf{u}_{n}$ from the set of transmit antennas at the time instant $n=1,\dots,N$. For simplicity, we still call this sequence of vectors as a sequence, which is denoted by $\mathbf{U}=\left[u_{n,m}\right]=\begin{bmatrix}\mathbf{u}_{1} & \cdots & \mathbf{u}_{N}\end{bmatrix}^{T}\in \mathbb{C}^{N\times N_{t}}$. Considering the unimodular constraint with energy budget $\mathrm{Tr}\left(\mathbf{U}^{H}\mathbf{U}\right)=\alpha$, we want to design $\mathbf{U}\in\mathcal{U}$, where
\begin{equation}
\mathcal{U}=\left\{ \mathbf{U}\middle| \left|u_{n,m}\right|= \sqrt{\frac{\alpha}{NN_{t}}},n=1,\ldots,N;m=1,\ldots,N_{t} \right\}.
\end{equation}
And the received sequence is given by
\begin{equation}
\label{eq:MIMO_channel}
	\mathbf{y}_{n}=\sum_{k=0}^{K}\mathbf{H}_{k}\mathbf{u}_{n-k}+\mathbf{v}_{n},
\end{equation}
where $\mathbf{u}_{n}=\mathbf{0}$ when $n\leq0$ or $n>N$, and $\mathbf{v}_{n}$ is an $N_{r}\times1$ noise vector. Equation \eqref{eq:MIMO_channel} can be written in a matrix form as
\begin{equation}
\label{eq:channel_model_matrix_form}
	\begin{bmatrix}\mathbf{y}_{1}^{T}\\ \vdots\\ \vdots\\ \vdots\\ \mathbf{y}_{N+K}^{T}\end{bmatrix}
	=
	\begin{bmatrix}
		\mathbf{u}_{1}^{T}
		&  & 
		\mathbf{0}\\ \vdots & \ddots & \mathbf{u}_{1}^{T}\\ \vdots & \ddots & \vdots\\
		\mathbf{u}_{N}^{T} & \ddots & \vdots\\ \mathbf{0}
		&  & 
		\mathbf{u}_{N}^{T}
	\end{bmatrix}
	\begin{bmatrix}\mathbf{H}_{0}^{T}\\ \vdots\\ \mathbf{H}_{K}^{T}\end{bmatrix}
	+
	\begin{bmatrix}\mathbf{v}_{1}^{T}\\ \vdots\\ \vdots\\ \vdots\\ \mathbf{v}_{N+K}^{T}\end{bmatrix}.
\end{equation}
Let $\mathbf{Y}=\begin{bmatrix}\mathbf{y}_{1} & \cdots & \mathbf{y}_{N+K}\end{bmatrix}^{T}\in\mathbb{C}^{(N+K)\times N_{r}}$ be the received matrix, and
\begin{equation}
\mathbf{S}=\mathcal{T}\left(\mathbf{U}\right)=
\begin{bmatrix}
\mathbf{u}_{1}^{T}
&  & 
\mathbf{0}\\ \vdots & \ddots & \mathbf{u}_{1}^{T}\\ \vdots & \ddots & \vdots\\
\mathbf{u}_{N}^{T} & \ddots & \vdots\\ \mathbf{0}
&  & 
\mathbf{u}_{N}^{T}
\end{bmatrix}
\in\mathbb{C}^{(N+K)\times \left(K+1\right)N_t}
\end{equation}
be a block Toeplitz convolution matrix with $\begin{bmatrix}\mathbf{U}^T & \mathbf{0}\end{bmatrix}^{T}$ being the first block and remaining blocks are obtained by a downward circular shift of the previous block. Note that since $\mathrm{Tr}\left(\mathbf{U}^{H}\mathbf{U}\right)=\alpha$, then $\mathrm{Tr}\left(\mathbf{S}^{H}\mathbf{S}\right)=\alpha (K+1)$.  $\mathbf{H}=\begin{bmatrix}\mathbf{H}_{0} & \cdots & \mathbf{H}_{K}\end{bmatrix}^{T}\in \mathbb{C}^{\left(K+1\right)N_{t}\times N_{r}}$
is the channel impulse response with matrix-form taps, and $\mathbf{V}$ is the noise matrix. Thus, we can write in a compact way the received signal as
\begin{equation}\label{eq:MIMO_channel_matform}
	\mathbf{Y}=\mathbf{S}\mathbf{H}+\mathbf{V}.
\end{equation}

It can be easily seen that each column of $\mathbf{Y}$ corresponds to a received sequence for one of the $N_{r}$ receive antennas, i.e., a multi-input single-output (MISO) channel. Let $\mathbf{y}=\mbox{vec}\left(\mathbf{Y}\right)$
, $\mathbf{h}=\mbox{vec}\left(\mathbf{H}\right)$, and $\mathbf{v}=\mbox{vec}\left(\mathbf{V}\right)$, and based on $\mathrm{vec}\left(\mathbf{XYZ}\right)=\left(\mathbf{Z}^T\otimes \mathbf{X}\right)\mathrm{vec}\left(\mathbf{Y}\right)$, we have
\begin{equation}
\label{eq:channel_model}
	\mathbf{y}=\left(\mathbf{I}_{N_{r}}\otimes\mathbf{S}\right)\mathbf{h}+\mathbf{v}.
\end{equation}

\subsection{Heuristic Existing Methods}
\label{ssec:Heuristic_Existing_Methods}

Most of the current works on unimodular sequence design focus on good auto- and cross-correlation properties; see \cite{he2009designing} on MIMO radar unimodular codes and references therein. The good correlation properties are particularly desired in that the matched filter is employed in subsequent channel estimation. As a matter of fact, the obtained channel estimate is closely related to maximum likelihood (ML) estimation. Assume the vectorized noise in the channel model \eqref{eq:channel_model} follows a circularly complex Gaussian distribution, $\mathbf{v}\sim\mathcal{CN}\left(\mathbf{0},\sigma^2\mathbf{I}\right)$. Minimizing the mean square error (MSE) $\mathrm{E}\{\|\hat{\mathbf{h}}_{\mathrm{ML}}-\mathbf{h}\|^{2}\}$ results in the ML channel estimate \cite{kay1993fundamentals}
\begin{align}
	\hat{\mathbf{h}}_{\mathrm{ML}}
	&= \left(\left(\mathbf{I}_{N_{r}}\otimes\mathbf{S}\right)^{H}\left(\mathbf{I}_{N_{r}}\otimes\mathbf{S}\right)\right)^{-1}\left(\mathbf{I}_{N_{r}}\otimes\mathbf{S}\right)^{H}\mathbf{y}\\
	&= \left(\mathbf{I}_{N_{r}}\otimes\mathbf{S}^{H}\mathbf{S}\right)^{-1}\left(\mathbf{I}_{N_{r}}\otimes\mathbf{S}\right)^{H}\mathbf{y},
\end{align}
where the second equality is due to $\left(\mathbf{X\otimes \mathbf{Y}}\right)\left(\mathbf{M}\otimes \mathbf{N}\right)=\mathbf{XM}\otimes\mathbf{YN}$. And the corresponding error is given by
\begin{align}
	\mathcal{E}
	&= \mathrm{Tr}\left(\left(\left(\mathbf{I}_{N_{r}}\otimes\mathbf{S}\right)^{H}\left(\mathbf{I}_{N_{r}}\otimes\mathbf{S}\right)\right)^{-1}\right)\\
	&=
	\mathrm{Tr}\left(\left(\mathbf{I}_{N_{r}}\otimes\mathbf{S}^{H}\mathbf{S}\right)^{-1}\right)\\
	&=
	N_{r}\mathrm{Tr}\left(\left(\mathbf{S}^{H}\mathbf{S}\right)^{-1}\right).\label{eq:ML_min_problem}
\end{align}
To minimize the error of ML estimation, the training sequence should be a solution to the optimization problem
\begin{equation}
\label{eq:ML_problem}
	\begin{array}{rcrc}
		\underset{\mathbf{U},\mathbf{S}}{\text{minimize}} & \mathcal{E} & \text{subject to} & \mathbf{S}=\mathcal{T}\left(\mathbf{U}\right),\mathbf{U}\in\mathcal{U}.
	\end{array}
\end{equation}
\begin{mylemma}[\cite{manton2001optimal}]\label{lemma1}
	Let $\mathbf{X}\in\mathbb{C}^{M\times N}$ be such that $\mathrm{Tr}\left(\mathbf{X}^{H}\mathbf{X}\right)\leq\mu$ for some constant $\mu$. The minimum of $\mathrm{Tr}\left(\left(\mathbf{X}^{H}\mathbf{X}\right)^{-1}\right)$ is achieved when $\mathbf{X}^{H}\mathbf{X}=\frac{\mu}{N}\mathbf{I}$, provided that inverse of $\mathbf{X}^{H}\mathbf{X}$ exists.
\end{mylemma}

An approximation to problem \eqref{eq:ML_problem} is as follows. According to Lemma \ref{lemma1}, the objective function \eqref{eq:ML_min_problem} is minimized when  ($\mathrm{Tr\left(\mathbf{S}^{H}\mathbf{S}\right)}=\alpha(K+1)$)
\begin{equation}
\label{eq:ML_opt_S}
	\mathbf{S}^{H}\mathbf{S}=\frac{\alpha}{N_{t}}\mathbf{I}_{\left(K+1\right)N_{t}},
\end{equation}
if only the energy constraint is considered. Therefore, a heuristic approximation of the ML optimal sequence design could be formulated as
\begin{equation}
	\begin{array}{rl}
	\label{eq:ML_problem_2}
		\underset{\mathbf{U},\mathbf{S}}{\text{minimize}} & \left\Vert \mathbf{S}^{H}\mathbf{S}-\frac{\alpha}{N_{t}}\mathbf{I}\right\Vert _{F}^{2}\\
		\text{subject to} & \mathbf{S}=\mathcal{T}\left(\mathbf{U}\right),\mathbf{U}\in\mathcal{U}.
	\end{array}
\end{equation}

The optimal $\mathbf{S}$ satisfying~\eqref{eq:ML_opt_S} portrays an impulse-like correlation shape pursued in \cite{he2009designing,song2015seqencce}, where the aperiodic cross-correlation is defined as
\begin{equation}
\label{eq:cross_corr}
	r_{m_{1},m_{2}}(k)=\sum_{n=k+1}^{N}u_{n,m_{1}}u_{n-k,m_{2}}^{\ast}
\end{equation}
for $m_{1},m_{2}=1,\dots,N_{t}$ and lags $k=0,\dots,N-1$. Equation \eqref{eq:cross_corr} also defines the autocorrelation for the sequence of each transmit antenna when $m_{1}=m_{2}$. Accordingly, the correlation matrices for different lags $k=-(N-1),\dots,0,\dots,(N-1)$ are given by
\begin{equation}
\label{eq:corr_at_lag}
	\boldsymbol{\Sigma}_{k}=\begin{bmatrix}r_{1,1}\left(k\right) & r_{1,2}\left(k\right) & \cdots & r_{1,N_{t}}\left(k\right)\\
	r_{2,1}\left(k\right) & r_{2,2}\left(k\right) & \cdots & r_{2,N_{t}}\left(k\right)\\
	\vdots & \vdots & \ddots & \vdots\\
	r_{N_{t},1}\left(k\right) & r_{N_{t},2}\left(k\right) & \cdots & r_{N_{t},N_{t}}\left(k\right)
	\end{bmatrix},
\end{equation}
with $r_{m_{1},m_{2}}(-k)=r_{m_{1},m_{2}}^\ast(k)$, and $\boldsymbol{\Sigma}_{-k}=\boldsymbol{\Sigma}_{k}^{H}$. Let us define the correlation matrix for a sequence $\mathbf{S}$ as
\begin{equation}
	\boldsymbol{\Sigma}=\mathbf{S}^{H}\mathbf{S},
\end{equation}
and then we have
\begin{equation}
	\boldsymbol{\Sigma}=\begin{bmatrix}\boldsymbol{\Sigma}_{0} & \boldsymbol{\Sigma}_{-1} & \cdots & \boldsymbol{\Sigma}_{-K}\\
	\boldsymbol{\Sigma}_{1} & \boldsymbol{\Sigma}_{0} & \cdots & \boldsymbol{\Sigma}_{-\left(K-1\right)}\\
	\vdots & \vdots & \ddots & \vdots\\
	\boldsymbol{\Sigma}_{K} & \boldsymbol{\Sigma}_{K-1} & \cdots & \boldsymbol{\Sigma}_{0}
	\end{bmatrix}.
\end{equation}
Note that $\boldsymbol{\Sigma}$ only describes correlations at lags of interest, which in this case is determined by the length of the channel impulse response. To achieve the optimality dictated by \eqref{eq:ML_opt_S}, we can rewrite approximation problem \eqref{eq:ML_problem_2} as
\begin{equation}
\label{eq:ML_problem_3}
	\begin{array}{rl}
		\underset{\mathbf{U}}{\text{minimize}} & (K+1)\left\Vert \boldsymbol{\Sigma}_{0}- \frac{\alpha}{N_{t}}\mathbf{I}\right\Vert _{F}^{2}\\
		& +2\sum\limits_{k=1}^{K}(K+1-k)\left\Vert \boldsymbol{\Sigma}_{k}\right\Vert _{F}^{2}\\
		\text{subject to} & \mathbf{U}\in\mathcal{U}.
	\end{array}
\end{equation}
The objective function of \eqref{eq:ML_problem_3} is indeed the weighted correlation minimization criterion within the lag interval $k=0,\dots,K$ \cite{he2009designing}, for which algorithms WeCan and CAD were proposed. Another formulation is also presented in a similar attempt to procure the good correlation property as
\begin{equation}
	\begin{array}{rl}
	\label{eq:ML_problem_4}
		\underset{\mathbf{U}}{\text{minimize}} & \left\Vert \boldsymbol{\Sigma}_{0}- \frac{\alpha}{N_{t}}\mathbf{I}\right\Vert _{F}^{2}+2\sum\limits_{k=1}^{N-1}\left\Vert \boldsymbol{\Sigma}_{k}\right\Vert _{F}^{2}\\
		\text{subject to} & \mathbf{U}\in\mathcal{U},
	\end{array}
\end{equation}
for which an algorithm called CAN was developed in \cite{stoica2009new}, and WeCAN can be employed as well. In \cite{song2015seqencce}, both problems \eqref{eq:ML_problem_3} and \eqref{eq:ML_problem_4} were studied by considering a more general weighted formulation, and efficient algorithms were proposed.

We can see that sequences with good auto- and cross-correlation properties are desirable in general as no prior information on the channel is taken into account in the ensuing channel estimation task. Channel statistics, however, are often available on both the transmitter sides and receiver sides, and incorporating those priors into the design of the training sequence will improve the performance of channel estimator. In the following subsections, we will formulate the unimodular sequence design problem based on the MMSE minimization and the CMI maximization, both of which have been adopted as criteria in various estimation problems. In order for the channel model~\eqref{eq:channel_model} to be general, we assume $\mathbf{h}\sim\mathcal{CN}\left(\mathbf{h}_{0},\mathbf{R}_{0}\right)$, and the noise $\mathbf{v}\sim\mathcal{CN}\left(\mathbf{0},\mathbf{W}\right)$, where both the channel covariance $\mathbf{R}_{0}$ and the noise covariance $\mathbf{W}$ are arbitrary.

\subsection{Optimal Sequence Design by Minimizing the MMSE Criterion}
Given the channel model \eqref{eq:channel_model}, by minimizing the MSE $\mathrm{E}\{\|\hat{\mathbf{h}}_{\mathrm{MMSE}}-\mathbf{h}\|^{2}\}$, the MMSE estimator of the channel impulse $\mathbf{h}$ is given by
\begin{equation}
\label{eq:MMSE_estimator}
	\hat{\mathbf{h}}_{\mathrm{MMSE}}=\mathbf{R}_{0}\tilde{\mathbf{S}}^{H}\left(\tilde{\mathbf{S}}\mathbf{R}_{0}\tilde{\mathbf{S}}^{H}+\mathbf{W}\right)^{-1}\left(\mathbf{y}-\tilde{\mathbf{S}}\mathbf{h}_{0}\right)+\mathbf{h}_{0},
\end{equation}
where $\tilde{\mathbf{S}}=\mathbf{I}_{N_{r}}\otimes\mathbf{S}$ \cite{kay1993fundamentals}. And the error covariance matrix is
\begin{align}
	\mathbf{R}
	&=
	\mathrm{E}\left\{\left(\hat{\mathbf{h}}_{\mathrm{MMSE}}-\mathbf{h}\right)\left(\hat{\mathbf{h}}_{\mathrm{MMSE}}-\mathbf{h}\right)^{H}\right\}\\
	&= \mathbf{R}_{0}-\mathbf{R}_{0}\tilde{\mathbf{S}}^{H}\left(\tilde{\mathbf{S}}\mathbf{R}_{0}\tilde{\mathbf{S}}^{H}+\mathbf{W}\right)^{-1}\tilde{\mathbf{S}}\mathbf{R}_{0}\label{eq:cov_matrix}\\
	&=
	\left(\mathbf{R}_{0}^{-1}+\tilde{\mathbf{S}}^{H}\mathbf{W}^{-1}\tilde{\mathbf{S}}\right)^{-1},
\end{align}
where the last equality is due to the matrix inversion lemma \cite{palomar2003joint}. The MMSE is thus given by
\begin{equation}
	\mathrm{MMSE}\left(\mathbf{S}\right)=\mathrm{Tr}\left(\mathbf{R}\right),
\end{equation}
and the following problem can be formulated
\begin{equation}
\label{eq:MMSE_prob}
	\begin{array}{rcrc}
		\underset{\mathbf{U},\mathbf{S}}{\text{minimize}} & \mathrm{MMSE}\left(\mathbf{S}\right) & \text{subject to} & \mathbf{S}=\mathcal{T}\left(\mathbf{U}\right),\mathbf{U}\in\mathcal{U},
	\end{array}
\end{equation}
which gives the optimal unimodular training sequence for the MMSE channel estimation.

\subsection{Optimal Sequence Design by Maximizing the CMI Criterion}
Apart from the MMSE criterion, another popular statistical measure in channel estimation is the conditional mutual information (CMI) between the channel impulse response and the received sequence, e.g., \cite{Kotecha2004}. The CMI is defined as
\begin{align}
\label{eq:CMI_1}
	\mathrm{CMI}\left(\mathbf{S}\right)
	& = I\left(\mathbf{h};\mathbf{y}\left|\mathbf{S}\right.\right)\\
	& = H\left(\mathbf{h}\right)-H\left(\mathbf{h}\left|\mathbf{y},\mathbf{S}\right.\right),
\end{align}
where $H\left(\cdot\right)$ is the differential entropy of a distribution \cite{cover2012elements}. Under the linear model \eqref{eq:channel_model} with Gaussian distributed channel impulse and noise, we have the conditional distribution $\mathbf{h}\left|\mathbf{y},\mathbf{S}\right.\sim\mathcal{CN}(\hat{\mathbf{h}},\mathbf{R})$. Then $\mathrm{CMI}\left(\mathbf{S}\right)$ can be written as
\setlength{\arraycolsep}{0.0em}
\begin{eqnarray}
	\mathrm{CMI}\left(\mathbf{S}\right)
	&{}={}& \frac{1}{2}\log\left(\left(2\pi e\right)^{(K+1)N_{t}N_{r}}\det\left(\mathbf{R}_{0}\right)\right)\nonumber\\
	&&{-}\:\frac{1}{2}\log\left(\left(2\pi e\right)^{(K+1)N_{t}N_{r}}\det\left(\mathbf{R}\right)\right)\\
	&{}={}&
	\frac{1}{2}\log\det\left(\mathbf{R}_{0}\mathbf{R}^{-1}\right).
\end{eqnarray}
\setlength{\arraycolsep}{5pt}
By maximizing $\mathrm{CMI}\left(\mathbf{S}\right)$ we reach the following optimization problem
\begin{equation}
	\begin{array}[t]{cc}
	\label{eq:CMI_prob}
		\underset{\mathbf{U},\mathbf{S}}{\text{maximize}} & \mathrm{CMI}\left(\mathbf{S}\right)\end{array}\begin{array}[t]{cc}
		\text{subject to} & \mathbf{S}=\mathcal{T}\left(\mathbf{U}\right),\mathbf{U}\in\mathcal{U}.
	\end{array}
\end{equation}

\textit{Remark}: It should be mentioned that channel model \eqref{eq:channel_model} includes the SISO channel as a special case. A lot of efforts have been made to construct unimodular sequences for SISO channels via either analytical methods or computational approaches. Apart from early works on binary sequences and polyphase sequences, e.g., \cite{golay1972class,zhang1993polyphase}, numerical algorithms are provided to design unimodular sequences of good correlation properties \cite{stoica2009new,song2014optimal}. Let $N_{t}=N_{r}=1$ and $\mathbf{u}$ denote the training sequence, then $\mathbf{S}=\mathcal{T}\left(\mathbf{u}\right)$ is a Toeplitz convolution matrix and expression \eqref{eq:corr_at_lag} reduces to a scalar that gives autocorrelations at different lags for the sequence $\mathbf{u}$. And similar formulations as \eqref{eq:ML_problem_3} and \eqref{eq:ML_problem_4} are proposed in order to obtain sequences of good autocorrelation properties. However, as we have seen in the previous discussion, the resulting channel estimate cannot benefit from the available knowledge of channel statistics. Therefore, designing optimal training sequences by minimizing the MMSE criterion or maximizing the CMI criterion will be beneficial in terms of final estimation performances. Without any modifications, formulations \eqref{eq:MMSE_prob} and \eqref{eq:CMI_prob} can be deployed in the context of SISO channels.

\section{Algorithms for Unimodular Sequence Design}
\label{sec:Algorithms_for_Optimal}
In this section, we develop efficient algorithms to solve problems \eqref{eq:MMSE_prob} and \eqref{eq:CMI_prob}. There is an extensive literature dealing with optimization problems of similar objective functions with only power constraint on $\mathbf{S}$ where, assuming some special structure for the prior covariance matrices of channel and noise, the problems are reformulated as power allocation with waterfilling-like solutions. In our formulations, however, it is not only the Toeplitz structure of $\mathbf{S}$ but also the tough unimodular constraint that prevents us from adopting the same approach.

It is worth mentioning that a possible approach to problem \eqref{eq:MMSE_prob} is a two-stage procedure \cite{yang2010mimo} related to correlation shaping. If the channel noise is independent and identically distributed, i.e., $\mathbf{W}=\sigma^{2}\mathbf{I}$ for some power density $\sigma^{2}$, the objective function becomes
\begin{align}
	\mathrm{MMSE}\left(\mathbf{S}\right)
	&=
	\mathrm{Tr}\left(\left(\mathbf{R}_{0}^{-1}+\frac{1}{\sigma^{2}}\tilde{\mathbf{S}}^{H}\tilde{\mathbf{S}}\right)^{-1}\right)\\
	&=
	\mathrm{Tr}\left(\left(\mathbf{R}_{0}^{-1}+\frac{1}{\sigma^{2}}\mathbf{I}_{N_{r}}\otimes \boldsymbol{\Sigma}\right)^{-1}\right),\label{eq:MMSE_iid_prob}
\end{align}
where the second equality follows from substitution of the correlation matrix $\boldsymbol{\Sigma}=\mathbf{S}^{H}\mathbf{S}$. Consider only the constraint $\mathrm{Tr}\left(\boldsymbol{\Sigma}\right)=(K+1)\alpha$ induced by the energy budget in \eqref{eq:MMSE_prob}, minimizing \eqref{eq:MMSE_iid_prob} with respect to $\boldsymbol{\Sigma}$ (instead of $\mathbf{S}$) can be rewritten as an SDP by resorting to the Schur-complement theorem \cite{boyd2004convex}, which yields the optimal correlation matrix $\boldsymbol{\Sigma}^{\star}$. Once $\boldsymbol{\Sigma}^{\star}$ is obtained, the problem boils down to recovering sequences from its correlation matrix, which is to solve the following approximation problem
\begin{equation}
	\begin{array}{rl}
	\label{eq:design_by_corr_problem}
		\underset{\mathbf{U},\mathbf{S}}{\text{minimize}} & \Vert\mathbf{S}^{H}\mathbf{S} - \boldsymbol{\Sigma}^{\star}\Vert_{F}\\
		\text{subject to} & \mathbf{S}=\mathcal{T}\left(\mathbf{U}\right),\mathbf{U}\in\mathcal{U},
	\end{array}
\end{equation}
if zero error is achievable. For a single sequence and without the unimodular constraint, problem \eqref{eq:design_by_corr_problem} can be tackled by means of filter design \cite{davidson2010enriching}. However, constructing a unimodular sequence that presents a prescribed correlation shape is challenging. As a special case, \cite{he2009designing} and \cite{song2015seqencce} have studied this problem only when the correlation matrix is an identity. On the other hand, it is not guaranteed that the objective in \eqref{eq:design_by_corr_problem} can reach zero when minimized. For example, when the number of sequences is relatively large for the training length, it is impossible to design sequences such that correlation matrix is an identity, i.e., auto- and cross-correlation cannot be made small simultaneously \cite{he2009designing}.

Therefore, it is advisable to solve problems \eqref{eq:MMSE_prob} and \eqref{eq:CMI_prob} directly with the colored noise considered. In the following, we will devise algorithms for both problems based on the majorization-minimization framework.

\subsection{Majorization-Minimization Framework}
The majorization-minimization, or MM method is a general framework for solving an optimization problem indirectly. In this section, we will briefly introduce the idea of the MM method for a minimization problem, and the details can be found in \cite{hunter2004tutorial,beck2009gradient}.

The MM method tackles a difficult optimization problem by solving a series of simple approximation problems. Given a minimization problem
\begin{eqnarray}
\label{eq:MM_original_prob}
	\begin{array}{rcrc}
		\underset{\mathbf{x}}{\text{minimize}} & f(\mathbf{x}) & \text{subject to} & \mathbf{x}\in\mathcal{X},
	\end{array}
\end{eqnarray}
and a feasible starting point $\mathbf{x}^{(0)}\in\mathcal{X}$, the MM method minimizes a sequence of surrogate functions $g\left(\mathbf{x},\mathbf{x}^{(t)}\right),t=0,1,\dots$ instead. Each surrogate function is a majorization function of $f(\mathbf{x})$ at $\mathbf{x}^{(t)}$ that satisfies:
\begin{align}
	g\left(\mathbf{x}^{(t)},\mathbf{x}^{(t)}\right) & = f\left(\mathbf{x}^{(t)}\right),\label{eq:MM_majorizaiton_1}\\
	g\left(\mathbf{x},\mathbf{x}^{(t)}\right) & \geq f\left(\mathbf{x}\right) \text{ for every }\mathbf{x}\in\mathcal{X},\label{eq:MM_majorizaiton_2}
\end{align}
and
\begin{equation}
\label{eq:MM_majorized_prob}
	\mathbf{x}^{(t+1)}\in \underset{\mathbf{x}\in\mathcal{X}}{\arg\min}\enspace g\left(\mathbf{x},\mathbf{x}^{(t)}\right).
\end{equation}
According to the rules \eqref{eq:MM_majorizaiton_1} and \eqref{eq:MM_majorizaiton_2}, we have
\begin{equation}\label{eq:MM_inequs}
	f\left(\mathbf{x}^{(t+1)}\right)\leq g\left(\mathbf{x}^{(t+1)},\mathbf{x}^{(t)}\right)\leq g\left(\mathbf{x}^{(t)},\mathbf{x}^{(t)}\right)=f\left(\mathbf{x}^{(t)}\right),
\end{equation}
and consequently, the MM method produces a sequence of points $\mathbf{x}^{(t)}$, for which the original objective function of \eqref{eq:MM_original_prob} is monotonically decreased. Provided that the objective function is bounded below, it is guaranteed that the MM algorithm will converge to a stationary point.

The key question is then how to find a good majorization function $g\left(\mathbf{x},\mathbf{x}^{(t)}\right)$ such that the resulting problems \eqref{eq:MM_majorized_prob} are easy to solve. Although there is no universal rule to determine the function $g\left(\mathbf{x},\mathbf{x}^{(t)}\right)$, the structure of the problem at hand can nevertheless provide helpful hints and some tricks are suggested in \cite{hunter2004tutorial}.

\subsection{MM-Based Algorithms}
\label{ssec:Algorithms_Opt_UniSeq}
Let us introduce $\mathbf{P}=\tilde{\mathbf{S}}\mathbf{R}_{0}\tilde{\mathbf{S}}^{H}+\mathbf{W}$, and by \eqref{eq:cov_matrix} the objective function for the MMSE minimization problem \eqref{eq:MMSE_prob} can be written as
\begin{equation}
\label{eq:MMSE_function}
	\mathrm{MMSE}\left(\mathbf{S}\right) =  \mathrm{Tr}\left(\mathbf{R}_{0}-\mathbf{R}_{0}\tilde{\mathbf{S}}^{H}\mathbf{P}^{-1}\tilde{\mathbf{S}}\mathbf{R}_{0}\right).
\end{equation}
\begin{mylemma}\label{lemma2}
	The function $f(\mathbf{X,Z})=\mathrm{Tr}\left(\mathbf{X}^{H}\mathbf{Z}^{-1}\mathbf{X}\right)$ is a matrix fractional function and is jointly convex in $\mathbf{Z}\succ 0$ and $\mathbf{X}$ \cite{boyd2004convex}.
\end{mylemma}

By Lemma \ref{lemma2}, $\mathrm{MMSE}\left(\mathbf{S}\right)=\mathrm{Tr}\left(\mathbf{R}_{0}-\mathbf{R}_{0}\tilde{\mathbf{S}}^{H}\mathbf{P}^{-1}\tilde{\mathbf{S}}\mathbf{R}_{0}\right)$ is jointly concave in $\{\tilde{\mathbf{S}},\mathbf{P}\}$ (recall that $\tilde{\mathbf{S}}=\mathbf{I}_{N_{r}}\otimes \mathbf{S}$). Since a concave function is upper-bounded by its supporting hyperplane, $\mathrm{MMSE}\left(\mathbf{S}\right)$ can be majorized as follows:
\setlength{\arraycolsep}{0.0em}%
\begin{eqnarray}
\hspace{-0.5cm}
	\mathrm{MMSE}\left(\mathbf{S}\right)
	&{}\leq{}&
	g_{\mathrm{MMSE}}\left(\mathbf{S},\mathbf{S}^{(t)}\right)\label{eq:MMSE_Majorization}\\
	&{}={}&
	\mathrm{MMSE}\left(\mathbf{S}^{(t)}\right){+} \mathrm{Tr}\bigg(\left(\mathbf{A}^{(t)}\right)^{H}\tilde{\mathbf{S}}\mathbf{R}_{0}\tilde{\mathbf{S}}^{H}\nonumber\\
	&  &\mathbf{A}^{(t)}\bigg)-2\mathrm{Re}\left\{ \mathrm{Tr}\left(\mathbf{R}_{0}\left(\mathbf{A}^{(t)}\right)^{H}\tilde{\mathbf{S}}\right)\right\},
\end{eqnarray}
\setlength{\arraycolsep}{5pt}%
where~$\tilde{\mathbf{S}}^{(t)}=\mathbf{I}_{N_{r}}\otimes\mathbf{S}^{(t)}$ with $\mathbf{S}^{(t)}=\mathcal{T}\left(\mathbf{U}^{(t)}\right)$, and $\mathbf{A}^{\left(t\right)}=\left(\tilde{\mathbf{S}}^{(t)}\mathbf{R}_{0}\left(\tilde{\mathbf{S}}^{(t)}\right)^{H}+\mathbf{W}\right)^{-1}\tilde{\mathbf{S}}^{(t)}\mathbf{R}_{0}$. To solve problem \eqref{eq:MMSE_prob}, it suffices to solve iteratively the following problem:
\begin{equation}
	\begin{array}{rl}
	\label{eq:MMSE_MM_1}
		\underset{\mathbf{U},\mathbf{S}}{\text{minimize}} & g_{\mathrm{MMSE}}\left(\mathbf{S},\mathbf{S}^{(t)}\right)\\
		\text{subject to} & \mathbf{S}=\mathcal{T}\left(\mathbf{U}\right),\mathbf{U}\in\mathcal{U},
	\end{array}
\end{equation}

For problem \eqref{eq:CMI_prob}, the objective function can be written as
\begin{equation}
\label{eq:CMI_function}
	\mathrm{CMI}\left(\mathbf{S}\right) = \frac{1}{2}\log\det\left(\mathbf{R}_{0}\left(\mathbf{R}_{0}-\mathbf{R}_{0}\tilde{\mathbf{S}}^{H}\mathbf{P}^{-1}\tilde{\mathbf{S}}\mathbf{R}_{0}\right)^{-1}\right).
\end{equation}
\begin{mylemma}
	\label{lemma3}
	Given a positive semidefinite matrix $\mathbf{M}$, the function $h(\mathbf{Z,X})=\mathbf{M}-\mathbf{M}\mathbf{X}^{H}\mathbf{Z}^{-1}\mathbf{X}\mathbf{M}$ is matrix concave over $\mathbf{X}$ of an appropriate size and $\mathbf{Z}\succ 0$ \cite{boyd2004convex}. Since $-\log\det\left(\cdot\right)$ is matrix convex and decreasing over positive definite cone, $-\log\det\left(\mathbf{M}-\mathbf{M}\mathbf{X}^{H}\mathbf{Z}^{-1}\mathbf{X}\mathbf{M}\right)$ is convex in $\{\mathbf{Z,X}\}$.
\end{mylemma}
Owing to Lemma \ref{lemma3}, $\mathrm{CMI}\left(\mathbf{S}\right)$ is jointly convex in $\{\tilde{\mathbf{S}},\mathbf{P}\}$, and we can obtain the following minorization
\setlength{\arraycolsep}{0.0em}%
\begin{eqnarray}
\label{eq:CMI_Majorization}
	\mathrm{CMI}\left(\mathbf{S}\right) & \geq & g_{\mathrm{CMI}}\left(\mathbf{S},\mathbf{S}^{(t)}\right)\\
	&{}={}& \mathrm{Re}\left\{ \mathrm{Tr}\left(\mathbf{R}_{0}\left(\mathbf{R}^{(t)}\right)^{-1}\left(\mathbf{A}^{(t)}\right)^{H}\tilde{\mathbf{S}}\right)\right\}\\
	&{}={}& -\frac{1}{2}\mathrm{Tr}\left(\left(\mathbf{R}^{(t)}\right)^{-1}\left(\mathbf{A}^{(t)}\right)^{H}\tilde{\mathbf{S}}\mathbf{R}_{0}\tilde{\mathbf{S}}^{H}\mathbf{A}^{(t)}\right)\nonumber\\
	& &{+}\:\mathrm{CMI}\left(\mathbf{S}^{(t)}\right),
\end{eqnarray}
\setlength{\arraycolsep}{5pt}%
where
\setlength{\arraycolsep}{0.0em}%
\begin{eqnarray}
	\mathbf{R}^{(t)}
	&{}={}& 
	\mathbf{R}_{0}-\mathbf{R}_{0}\left(\tilde{\mathbf{S}}^{(t)}\right)^{H}\left(\tilde{\mathbf{S}}^{(t)}\mathbf{R}_{0}\left(\tilde{\mathbf{S}}^{(t)}\right)^{H}+\mathbf{W}\right)^{-1}\nonumber\\
	&&\tilde{\mathbf{S}}^{(t)}\mathbf{R}_{0}\\
	&{}={}& \mathbf{R}_{0}^{-1}+\left(\tilde{\mathbf{S}}^{(t)}\right)^{H}\mathbf{W}^{-1}\tilde{\mathbf{S}}^{(t)}.
\end{eqnarray}
\setlength{\arraycolsep}{5pt}%
As a result, solving the CMI maximization problem \eqref{eq:CMI_prob} is equivalent to solving the series of minorized problems
\begin{equation}
\label{eq:CMI_MM_1}
\begin{array}{rl}
	\underset{\mathbf{U},\mathbf{S}}{\text{maximize}} & g_{\mathrm{CMI}}\left(\mathbf{S},\mathbf{S}^{(t)}\right)\\
	\text{subject to} & \mathbf{S}=\mathcal{T}\left(\mathbf{U}\right),\mathbf{U}\in\mathcal{U}.
\end{array}
\end{equation}

Notice that problems \eqref{eq:MMSE_MM_1} and \eqref{eq:CMI_MM_1} share a similar form of objective function. Let
\setlength{\arraycolsep}{0.0em}%
\begin{eqnarray}
\label{eq:MM_obj}
\hspace{-0.6cm}
	g\left(\mathbf{S};\mathbf{S}^{(t)},\mathbf{V}^{(t)}\right)
	&{}={}& \mathrm{Tr}\left(\mathbf{V}^{(t)}\left(\mathbf{A}^{(t)}\right)^{H}\tilde{\mathbf{S}}\mathbf{R}_{0}\tilde{\mathbf{S}}^{H}\mathbf{A}^{(t)}\right)\nonumber\\
	&&{-}\:2\mathrm{Re}\left\{ \mathrm{Tr}\left(\mathbf{R}_{0}\mathbf{V}^{(t)}\left(\mathbf{A}^{(t)}\right)^{H}\tilde{\mathbf{S}}\right)\right\},
\end{eqnarray}
\setlength{\arraycolsep}{5pt}%
where $\mathbf{V}^{(t)}=\mathbf{I}$
for the MMSE minimization problem and $\mathbf{V}^{\left(t\right)}=\left(\mathbf{R}^{\left(t\right)}\right)^{-1}$
for the CMI maximization problem. After reversing the sign of objective function of \eqref{eq:CMI_MM_1} and ignoring the constants and the scaling factor, the following unified problem is obtained
\begin{equation}
	\begin{array}{rl}
	\label{eq:MM_prob}
		\underset{\mathbf{U},\mathbf{S}}{\text{minimize}} & g\left(\mathbf{S};\mathbf{S}^{(t)},\mathbf{V}^{(t)}\right)\\
		\text{subject to} & \mathbf{S}=\mathcal{T}\left(\mathbf{U}\right),\mathbf{U}\in\mathcal{U}.
	\end{array}
\end{equation}

\begin{mylemma}\label{lemma4}
	Given Hermitian $\mathbf{M}\in\mathbb{C}^{n\times n}$ and $\mathbf{Z}\in\mathbb{C}^{m\times m}$ and any $\mathbf{X}^{(t)}\in\mathbb{C}^{m\times n}$, the function $\mathrm{Tr}\left(\mathbf{Z}\mathbf{X}\mathbf{M}\mathbf{X}^{H}\right)$ can be majorized by $-2\mathrm{Re}\left\{ \mathrm{Tr}\left(\left(\lambda\mathbf{X}^{\left(t\right)}-\mathbf{Z}\mathbf{X}^{\left(t\right)}\mathbf{M}\right)^{H}\mathbf{X}\right)\right\}+\lambda\left\Vert \mathbf{X}\right\Vert _{F}^{2}+\mathrm{const}$, where $\lambda\mathbf{I}\succeq\mathbf{M}^{T}\otimes\mathbf{Z}$ for some constant $\lambda$.
	\begin{proof}
		Given $\lambda\mathbf{I}\succeq\mathbf{M}^{T}\otimes\mathbf{L}\mathbf{L}^{H}$ for some constant $\lambda$, we have
		\begin{align}
		\setlength{\arraycolsep}{0.0em}%
		& \mathrm{Tr}\left(\mathbf{Z}\mathbf{X}\mathbf{M}\mathbf{X}^{H}\right)\nonumber\\
		{}={}& \mathrm{vec}^{H}\left(\mathbf{X}\right)\mathrm{vec}\left(\mathbf{Z}\mathbf{X}\mathbf{M}\right)\\
		{}={}& \mathrm{vec}^{H}\left(\mathbf{X}\right)\left(\mathbf{M}^{T}\otimes\mathbf{Z}\right)\mathrm{vec}\left(\mathbf{X}\right)\\
		{}\leq{}&{-}2\mathrm{Re}\left\{ \mathrm{vec}^{H}\left(\mathbf{X}\right)\left(\lambda\mathbf{I}-\mathbf{M}^{T}\otimes\mathbf{Z}\right)\mathrm{vec}\left(\mathbf{X}^{\left(t\right)}\right)\right\}\nonumber\\
		& {+}\:\mathrm{vec}^{H}\left(\mathbf{X}^{(t)}\right)\left(\lambda\mathbf{I}-\mathbf{M}^{T}\otimes\mathbf{Z}\right)\mathrm{vec}\left(\mathbf{X}^{(t)}\right)\nonumber\\
		&{+}\:\lambda\mathrm{vec}^{H}\left(\mathbf{X}\right)\mathrm{vec}\left(\mathbf{X}\right)\\
		{}={}&{-}2\mathrm{Re}\left\{ \mathrm{Tr}\left(\lambda\mathbf{X}^{H}\mathbf{X}^{\left(t\right)}-\mathbf{Z}\mathbf{X}^{\left(t\right)}\mathbf{M}\mathbf{X}^{H}\right)\right\}\nonumber\\
		&{+}\:\lambda\left\Vert \mathbf{X}\right\Vert _{F}^{2}+\mathrm{vec}^{H}\left(\mathbf{X}^{(t)}\right)\left(\lambda\mathbf{I}-\mathbf{M}^{T}\otimes\mathbf{Z}\right)\mathrm{vec}\left(\mathbf{X}^{(t)}\right). 
		\setlength{\arraycolsep}{5pt}%
		\end{align}%
		Notice that the third term of the last equation is simply a constant. And a scalar version of Lemma \ref{lemma4} can be found in \cite[Lemma 1]{song2015seqencce}. 
	\end{proof}
\end{mylemma}

To solve problem \eqref{eq:MM_prob}, yet a second majorization can be applied with Lemma \ref{lemma4} (note that $\| \tilde{\mathbf{S}}\| ^2=N_{r}(K+1)\alpha$):
\setlength{\arraycolsep}{0.0em}%
\begin{eqnarray}
	&& g\left(\mathbf{S};\mathbf{S}^{(t)},\mathbf{V}^{(t)}\right)\nonumber\\
	&{}\leq{}
	&{-}2\mathrm{Re}\bigg\{\mathrm{Tr}\bigg(\lambda^{(t)}\tilde{\mathbf{S}}^{H}\tilde{\mathbf{S}}^{(t)} - \mathbf{A}^{(t)}\mathbf{V}^{(t)}\left(\mathbf{A}^{(t)}\right)^{H}\tilde{\mathbf{S}}^{(t)}\mathbf{R}_{0}\nonumber\\
	&&\tilde{\mathbf{S}}^{H}\bigg)\bigg\} - 2\mathrm{Re}\left\{ \mathrm{Tr}\left(\mathbf{R}_{0}\mathbf{V}^{(t)}\left(\mathbf{A}^{(t)}\right)^{H}\tilde{\mathbf{S}}\right)\right\} + \mathrm{const}\ \ \ \ \label{eq:Majorization_2}\\
	&{}={}
	&{-}2\mathrm{Re}\bigg\{\mathrm{Tr}\Big(\big(\lambda^{(t)}\tilde{\mathbf{S}}^{(t)} - \mathbf{A}^{(t)}\mathbf{V}^{(t)}\left(\mathbf{A}^{(t)}\right)^{H}\tilde{\mathbf{S}}^{(t)}\mathbf{R}_{0}\nonumber\\
	&&{+}\:\mathbf{A}^{(t)}\mathbf{V}^{(t)}\mathbf{R}_{0}\big)^{H}\tilde{\mathbf{S}}\Big)\bigg\} + \mathrm{const},
\end{eqnarray}
\setlength{\arraycolsep}{5pt}%
where $\lambda^{(t)}\mathbf{I}\succeq\mathbf{R}_{0}^{T}\otimes \mathbf{A}^{(t)}\mathbf{V}^{(t)}\left(\mathbf{A}^{(t)}\right)^{H}$. The tightest upper bound will be $\lambda^{(t)}=\lambda_{\mathrm{max}}\left(\mathbf{R}_{0}^{T}\otimes \mathbf{A}^{(t)}\mathbf{V}^{(t)}\left(\mathbf{A}^{(t)}\right)^{H}\right)$. But computing the largest eigenvalue is costly especially when the size of the matrix is large, and thus an alternative is advisable. Since both $\mathbf{R}_{0}$ and $\mathbf{A}^{(t)}\mathbf{V}^{(t)}\left(\mathbf{A}^{(t)}\right)^{H}$ are positive semidefinite matrices, the largest eigenvalues are bounded as
\setlength{\arraycolsep}{0.0em}%
\begin{align}
	\lambda_{\mathrm{max}}\left(\mathbf{R}_{0}\right)&{}\leq{}\left\|\mathbf{R}_{0}\right\|_{1},\\
	\hspace{-0.1cm}\lambda_{\mathrm{max}}\left(\mathbf{A}^{(t)}\mathbf{V}^{(t)}\left(\mathbf{A}^{(t)}\right)^{H}\right) &{}\leq{}\left\|\mathbf{A}^{(t)}\mathbf{V}^{(t)}\left(\mathbf{A}^{(t)}\right)^{H}\right\|_{1},
\end{align}
\setlength{\arraycolsep}{5pt}%
where $\|\cdot\|_{1}$ is maximum column sum matrix norm \cite{horn2012matrix}. With  $\lambda_{\mathrm{max}}\left(\mathbf{X\otimes\mathbf{Z}}\right)=\lambda_{\mathrm{max}}\left(\mathbf{X}\right)\lambda_{\mathrm{max}}\left(\mathbf{Z}\right)$, we propose
\begin{equation}
	\lambda^{(t)}=\left\|\mathbf{R}_{0}\right\|_{1}\left\|\mathbf{A}^{(t)}\mathbf{V}^{(t)}\left(\mathbf{A}^{(t)}\right)^{H}\right\|_{1}.
\end{equation}

Let $\mathbf{B}\left(\tilde{\mathbf{S}}^{(t)},\mathbf{V}^{(t)}\right)=\lambda^{(t)}\tilde{\mathbf{S}}^{(t)}-\mathbf{A}^{(t)}\mathbf{V}^{(t)}\left(\mathbf{A}^{(t)}\right)^{H}\tilde{\mathbf{S}}^{(t)}\mathbf{R}_{0}\\+\mathbf{A}^{(t)}\mathbf{V}^{(t)}\mathbf{R}_{0}$, and considering $\tilde{\mathbf{S}}=\mathbf{I}_{N_{r}}\otimes\mathbf{S}$ with $\mathbf{S}=\mathcal{T}\left(\mathbf{U}\right)$, we have
\begin{equation}
\label{eq:MM_prob_2}
	\mathbf{U}^{(t+1)}\in \underset{|u_{n,m}|=\sqrt{\frac{\alpha}{NN_{t}}}}{\arg\min}\enspace -2\mathrm{Re}\Bigg\{ \mathrm{Tr}\bigg(\Big(\sum_{i,j}\mathbf{B}[i,j]\Big)^{H}\mathbf{U}\bigg)\Bigg\},
\end{equation}
where $\mathbf{B}[i,j]$ is a submatrix of $\mathbf{B}$ with rows from $(N+K)(i-1)+j$ to $(N+K)(i-1)+N+j-1$ and columns from $N_{t}(K+1)(i-1)+N_{t}(j-1)+1$ to $N_{t}(K+1)(i-1)+N_{t}j$, for $i=1,\dots,N_{r}$, and $j=1,\dots,K+1$. To find the next update $\mathbf{U}^{(t+1)}$, note that \eqref{eq:MM_prob_2} can be equivalently written as
\begin{equation}
	\mathbf{U}^{(t+1)}\in \underset{|u_{n,m}|=\sqrt{\frac{\alpha}{NN_{t}}}}{\arg\min}\enspace \left\|\mathbf{U}-\sum\nolimits_{i,j}\mathbf{B}[i,j] \right\|_{F}^{2}.
\end{equation}
And the minimum is achieved by projection onto a complex circle, which is
\begin{equation}
	\mathbf{U}^{(t+1)}=\sqrt{\frac{\alpha}{NN_{t}}}e^{j\arg\left(\sum\nolimits_{i,j}\mathbf{B}[i,j]\right)},
\end{equation}
where $\arg(\cdot)$ is taken element-wise. The whole procedure is summarized in Algorithm \ref{alg:Optimal_unimodular_sequence}. The iterations of the algorithm is deemed to be converged, e.g., when the difference between two consecutive updates for $\mathbf{U}$ is no larger than some admitted threshold.
\begin{algorithm}[H]
	\caption{Design of unimodular training sequence for the MMSE minimization \eqref{eq:MMSE_prob} or the CMI maximization \eqref{eq:CMI_prob}.}
	\label{alg:Optimal_unimodular_sequence}
	\begin{algorithmic}[1]
		\State Set $t=0$,  and initialize $u^{(0)}_{n,m},n=1,\ldots,N;m=1,\ldots,N_{t}$.
		\Repeat 
		\State $\mathbf{S}^{(t)}=\mathcal{T}\left(\mathbf{U}^{(t)}\right)$, and $\tilde{\mathbf{S}}^{(t)}=\mathbf{I}_{N_{r}}\otimes\mathbf{S}^{(t)}$\label{eq:step3}
		\State $\mathbf{A}^{(t)}=\left(\tilde{\mathbf{S}}^{(t)}\mathbf{R}_{0}\left(\tilde{\mathbf{S}}^{(t)}\right)^{H}+\mathbf{W}\right)^{-1}\tilde{\mathbf{S}}^{(t)}\mathbf{R}_{0}$
		\State 
		$\mathbf{V}^{(t)}=
		\begin{cases}
		\mathbf{I}, & \text{for the MMSE minimization}\\
		\mathbf{R}^{(t)}, & \text{for the CMI maximization}
		\end{cases}$
		\State $\lambda^{(t)}=\left\|\mathbf{R}_{0}\right\|_{1}\left\|\mathbf{A}^{(t)}\mathbf{V}^{(t)}\left(\mathbf{A}^{(t)}\right)^{H}\right\|_{1}$ 
		\State \parbox[t]{\dimexpr\linewidth-\algorithmicindent}{$\mathbf{B}\left(\tilde{\mathbf{S}}^{(t)},\mathbf{V}^{(t)}\right)=\lambda^{(t)}\tilde{\mathbf{S}}^{(t)}-\mathbf{A}^{(t)}\mathbf{V}^{(t)}\left(\mathbf{A}^{(t)}\right)^{H}\tilde{\mathbf{S}}^{(t)}\mathbf{R}_{0}\\+\mathbf{A}^{(t)}\mathbf{V}^{(t)}\mathbf{R}_{0}$\strut}
		\State $\mathbf{U}^{(t+1)}=\sqrt{\frac{\alpha}{NN_{t}}}e^{j\arg\left(\sum\nolimits_{i,j}\mathbf{B}[i,j]\right)}$\label{eq:step8}
		\State $t\gets t+1$
		\Until convergence
	\end{algorithmic}
\end{algorithm}

\subsection{Convergence Analysis}
\label{ssec:convergence_analysis}
Algorithm \ref{alg:Optimal_unimodular_sequence} is essentially based on the majorization-minimization framework, which has been shown to converge to a stationary point for bounded objective functions. The generated sequence of points ${\mathbf{U}^{(t)}},t=0,1,\dots,$ monotonically decreases or increases the objective function for minimization and maximization problems, respectively. In this section, we give a detailed analysis of the convergence for Algorithm \ref{alg:Optimal_unimodular_sequence}. Without loss of generality, we only consider minimizing the MMSE criterion.

For a constrained minimization problem with a smooth objective function, a stationary point is obtained when the following first-order optimality condition is satisfied.
\begin{myproposition}\label{proposition_1}
	Let $f:\mathbb{R}^{N}\to \mathbb{R}$ be a smooth function. 
	A point $\mathbf{x}^{\star}$ is a local minimum of $f$ within a subset $\mathcal{X}\subset \mathbb{R}^{N}$ if
	\begin{equation}
		\nabla f(\mathbf{x}^{\star})^{T}\mathbf{y}\geq 0,\forall \mathbf{y}\in T_{\mathcal{X}}(\mathbf{x}^{\star}),
	\end{equation}
	where $T_{\mathcal{X}}(\mathbf{x}^{\star})$ is the tangent cone of $\mathcal{X}$ at $\mathbf{x}^{\star}$.
\end{myproposition}

Provided Proposition \ref{proposition_1}, the convergence of our proposed algorithm is guaranteed as follows.
\begin{mytheorem}\label{theorem_1}
	By solving the series of problems \eqref{eq:MM_prob_2} in Algorithm \ref{alg:Optimal_unimodular_sequence}, a sequence of points $\{\mathbf{U}^{(t)},t=0,\dots\}$ is obtained, of which every limit point is a stationary point of problem \eqref{eq:MMSE_prob}.
	\begin{proof}
		A similar proof has been given in \cite{song2014optimal}. For details please refer to \cite[Theorem 5]{song2014optimal}.
	\end{proof}
\end{mytheorem}

\subsection{Accelerated Algorithm}
\label{ssec:acceleration}
To develop Algorithm \ref{alg:Optimal_unimodular_sequence} for solving problems \eqref{eq:MMSE_prob} and \eqref{eq:CMI_prob}, the original function was majorized/minorized twice, which may result in a loose surrogate function; see \eqref{eq:MMSE_Majorization}, \eqref{eq:CMI_Majorization} and \eqref{eq:Majorization_2}. And the performance of the MM method is susceptible to the slow convergence as EM-like algorithms. Then following the same idea in \cite{song2014optimal,song2015seqencce}, we employ an off-the-shelf method, called squared iterative methods (SQUAREM) \cite{varadhan2008simple}, to accelerate Algorithm \ref{alg:Optimal_unimodular_sequence}. SQUAREM was originally proposed to improve the convergence of EM-type algorithms and simultaneously keep its simplicity and stability. It can be easily applied to accelerate the MM algorithms as well. For details of convergence analysis, also refer to \cite{varadhan2008simple}. Without loss of generality, we only consider acceleration of Algorithm \ref{alg:Optimal_unimodular_sequence} for the MMSE minimization problem. For the CMI maximization problem, a similar procedure can be followed. 

Given the current point $\mathbf{U}^{(t)}$, we call iterative steps \ref{eq:step3} to \ref{eq:step8} of Algorithm 1 collectively as one MM update, denoted by $\text{MMupdate}(\mathbf{U}^{(t)})$. The accelerated computing scheme is given by Algorithm \ref{alg:Accelerated_MM}. The step length is chosen by the Cauchy-Barzilai-Borwein (CBB) method. And the back-tracking step is adopted to maintain the monotone property of generated iterates. To guarantee its feasibility, projection to the constrained set $\mathcal{U}$ in steps 8 and 9 are applied.
\begin{algorithm}
	\caption{Accelerated scheme for designing optimal unimodular training sequence for the MMSE estimation.}
	\label{alg:Accelerated_MM}
	\begin{algorithmic}[1]
		\State Set $t=0$,  and initialize $u^{(0)}_{n,m},n=1,\ldots,N;m=1,\ldots,N_{t}$.
		\Repeat		
		\State $\mathbf{U}_{1}=\mbox{MMupdate}\left(\mathbf{U}^{(t)}\right)$
		\State $\mathbf{U}_{2}=\mbox{MMupdate}\left(\mathbf{U}_{1}\right)$
		\State $\mathbf{L}_{1}=\mathbf{U}_{1}-\mathbf{U}^{(t)}$
		\State $\mathbf{L}_{2}=\mathbf{U}_{2}-\mathbf{U}_{1}-\mathbf{L}_{1}$
		\State Step length $l=-\frac{\left\Vert \mathbf{L}_{1}\right\Vert _{F}}{\left\Vert \mathbf{L}_{2}\right\Vert _{F}}$
		\State $\mathbf{U}^{(t+1)}=\sqrt{\frac{\alpha}{NN_{t}}}e^{j\arg\left(\mathbf{U}^{(t)}-2l\mathbf{L}_{1}+l^{2}\mathbf{L}_{2}\right)}$\label{algorithm_2_step}
		\While{$\mathrm{MMSE}\left(\mathbf{S}^{(t+1)}\right)>\mathrm{MMSE}\left(\mathbf{S}^{(t)}\right)$}
		\State $l\gets\frac{l-1}{2}$, and go to step \ref{algorithm_2_step}
		\EndWhile
		\State $t\gets t+1$
		\Until{convergence} 
	\end{algorithmic}
\end{algorithm}

\section{Algorithms for Sequence Design under PAR Constraints}
\label{sec:Opt_Seq_PAR}
The unimodular constraint on the training sequence originates partly from the low peak-to-average power ratio (PAR) demand, e.g., in MIMO radar systems. Low PAR sequences have found many applications in practice because they can mitigate the non-linear effects at the transmitter side while enabling more flexibility of the designed sequences compared with unimodular ones. In this section, we consider the problem of designing optimal sequences with low PAR.

For a sequence of vectors $\mathbf{U}\in \mathbb{C}^{N\times N_{t}}$, $\mathbf{U}_{:,m}$ denotes the length-$N$ sequence sent from the $m$th antenna, for $m=1,\dots,N_{t}$. And PAR is usually defined for each sequence transmitted by a single antenna as
\begin{equation}
	\mathrm{PAR}(\mathbf{U}_{:,m})=\frac{\underset{n}{\mathrm{max}}\{|u_{n,m}|^{2}\}}{\frac{1}{N}\alpha_{m}},
\end{equation}
provided that the training energy for the $m$th antenna is $\|\mathbf{U}_{:,m}\|^{2} = \alpha_{m}$. Determining training energy for each transmit antenna may depend on power distribution among antennas satisfying $\sum_{m=1}^{N_{t}}\alpha_{m}=\alpha$. And it follows that $1\leq \mathrm{PAR}(\mathbf{U}_{:,m})\leq N$. When $\mathrm{PAR}(\mathbf{U}_{:,m})=1$, PAR constraint reduces to the unimodular constraint. Given the PAR constraints for each transmit antenna
\begin{equation}
	\mathrm{PAR}(\mathbf{U}_{:,m})\leq\xi_{m},m=1,\dots,N_{t}
\end{equation}
the optimal sequence design problem for minimizing MMSE is then formulated as
\begin{equation}
\label{eq:MMSE_prob_PAR}
	\begin{array}{rl}
		\underset{\mathbf{U},\mathbf{S}}{\text{minimize}} & \mathrm{MMSE}\left(\mathbf{S}\right)\\ \text{subject to} & \mathbf{S}=\mathcal{T}\left(\mathbf{U}\right)\\
		& \|\mathbf{U}_{:,m}\|^{2}=\alpha_{m}\\
		& \underset{n}{\text{max}}\{|u_{n,m}|\}\leq\sqrt{\frac{\alpha_{m}\xi_{m}}{N}},m=1,\dots,N_{t}
	\end{array}
\end{equation}
where $\mathrm{MMSE}(\mathbf{S})$ is given by \eqref{eq:MMSE_function}. For the CMI maximization, an optimization problem can be similarly formulated, which maximizes $\mathrm{CMI}(\mathbf{S})$ \eqref{eq:CMI_function} under the same constraints as that of \eqref{eq:MMSE_prob_PAR}.

Following the same procedure of applying the MM framework in Section \ref{ssec:Algorithms_Opt_UniSeq}, the following majorized (minorized) problems can be obtained for problem \eqref{eq:MMSE_prob_PAR} for the MMSE minimization (CMI maximization)
\begin{equation}
\begin{array}{rl}
\label{eq:MM_prob_PAR_1}
	\underset{\mathbf{U}}{\text{minimize}} & \left\|\mathbf{U}-\sum\nolimits_{i,j}\mathbf{B}_{i,j} \right\|_{F}^{2}\\
	\text{subject to} & \|\mathbf{U}_{:,m}\|^{2}=\alpha_{m}\\
	& \underset{n}{\text{max}}\{|u_{n,m}|\}\leq\sqrt{\frac{\alpha_{m}\xi_{m}}{N}},m=1,\dots,N_{t}\\
\end{array}
\end{equation}
It is obvious that problem \eqref{eq:MM_prob_PAR_1} can be separated into $N_{t}$ problems as
\begin{equation}
	\begin{array}{rl}
	\label{eq:MM_prob_PAR_2}
		\underset{\mathbf{U}_{:,m}}{\text{minimize}} & \left\|\mathbf{U}_{:,m}-\mathbf{c}_{m} \right\|^{2}\\
		\text{subject to} & \|\mathbf{U}_{:,m}\|^{2}=\alpha_{m}\\
		& \underset{n}{\text{max}}\{|u_{n,m}|\}\leq\sqrt{\frac{\alpha_{m}\xi_{m}}{N}},
	\end{array}
\end{equation}
for $m=1,\dots,N_{t}$, where $\mathbf{c}_{m}$ is the $m$th column of $\sum\nolimits_{i,j}\mathbf{B}_{i,j}$. Problem \eqref{eq:MM_prob_PAR_2} is a nearest vector problem with low PAR constraint and has been well studied in \cite{tropp2005designing} via Karush-Kuhn-Tucker (KKT) conditions. By using the well-developed algorithms in \cite{tropp2005designing} to solve each problem \eqref{eq:MM_prob_PAR_2}, the overall algorithm is summarized in Algorithm \ref{alg:Optimal_PAR_sequence}. Note that Algorithm \ref{alg:Optimal_PAR_sequence} shares the same convergence property as that of Algorithm \ref{alg:Optimal_unimodular_sequence}. Furthermore, the acceleration scheme based on the SQUAREM method is also applicable here, and the procedure is similar to Algorithm \ref{alg:Accelerated_MM}.

\begin{algorithm}[H]
\caption{Design of optimal training sequence for the MMSE minimization \eqref{eq:MMSE_prob} or the CMI maximization \eqref{eq:CMI_prob} under the PAR constraint.}
\label{alg:Optimal_PAR_sequence}
\begin{algorithmic}[1]
	\State Set $t=0$,  and initialize $\mathbf{U}^{(0)}$ such that $\underset{n}{\text{max}}\enspace \{|u^{(0)}_{n,m}|\}\leq\sqrt{\frac{\alpha_{m}}{N}},m=1,\dots,N_{t}$.
	\Repeat		
	\State $\mathbf{S}^{(t)}=\mathcal{T}\left(\mathbf{U}^{(t)}\right)$, and $\tilde{\mathbf{S}}^{(t)}=\mathbf{I}_{N_{r}}\otimes\mathbf{S}^{(t)}$
	\State $\mathbf{A}^{(t)}=\left(\tilde{\mathbf{S}}^{(t)}\mathbf{R}_{0}\left(\tilde{\mathbf{S}}^{(t)}\right)^{H}+\mathbf{W}\right)^{-1}\tilde{\mathbf{S}}^{(t)}\mathbf{R}_{0}$
	\State 
	$\mathbf{V}^{(t)}=
	\begin{cases}
	\mathbf{I}, & \text{for the MMSE minimization}\\
	\mathbf{R}^{(t)}, & \text{for the CMI maximization}
	\end{cases}$
	\State $\lambda^{(t)}=\left\|\mathbf{R}_{0}\right\|_{1}\left\|\mathbf{A}^{(t)}\mathbf{V}^{(t)}\left(\mathbf{A}^{(t)}\right)^{H}\right\|_{1}$ 
	\State \parbox[t]{\dimexpr\linewidth-\algorithmicindent}{$\mathbf{B}\left(\tilde{\mathbf{S}}^{(t)},\mathbf{V}^{(t)}\right)=\lambda^{(t)}\tilde{\mathbf{S}}^{(t)}-\mathbf{A}^{(t)}\mathbf{V}^{(t)}\left(\mathbf{A}^{(t)}\right)^{H}\tilde{\mathbf{S}}^{(t)}\mathbf{R}_{0}\\+\mathbf{A}^{(t)}\mathbf{V}^{(t)}\mathbf{R}_{0}$\strut} 
	\State \parbox[t]{\dimexpr\linewidth-\algorithmicindent}
	{$\mathbf{U}^{(t+1)}_{:,m}\in \underset{\begin{subarray}{c}{\text{max}}_{n}\{|u_{n,m}|\}\leq\sqrt{\frac{\alpha_{m}\xi_{m}}{N}}\\ \|\mathbf{U}_{:,m}\|^{2}=\alpha_{m}\end{subarray}}{\arg\min}\enspace \left\|\mathbf{U}_{:,m}-\mathbf{c}_{m}\right\|^{2},m=1,\dots,N_{t}$\strut}
	\State $t\gets t+1$
	\Until{convergence} 
\end{algorithmic}
\end{algorithm}

\begin{figure}[!t]
	\centering
	\includegraphics[width=3.4in]{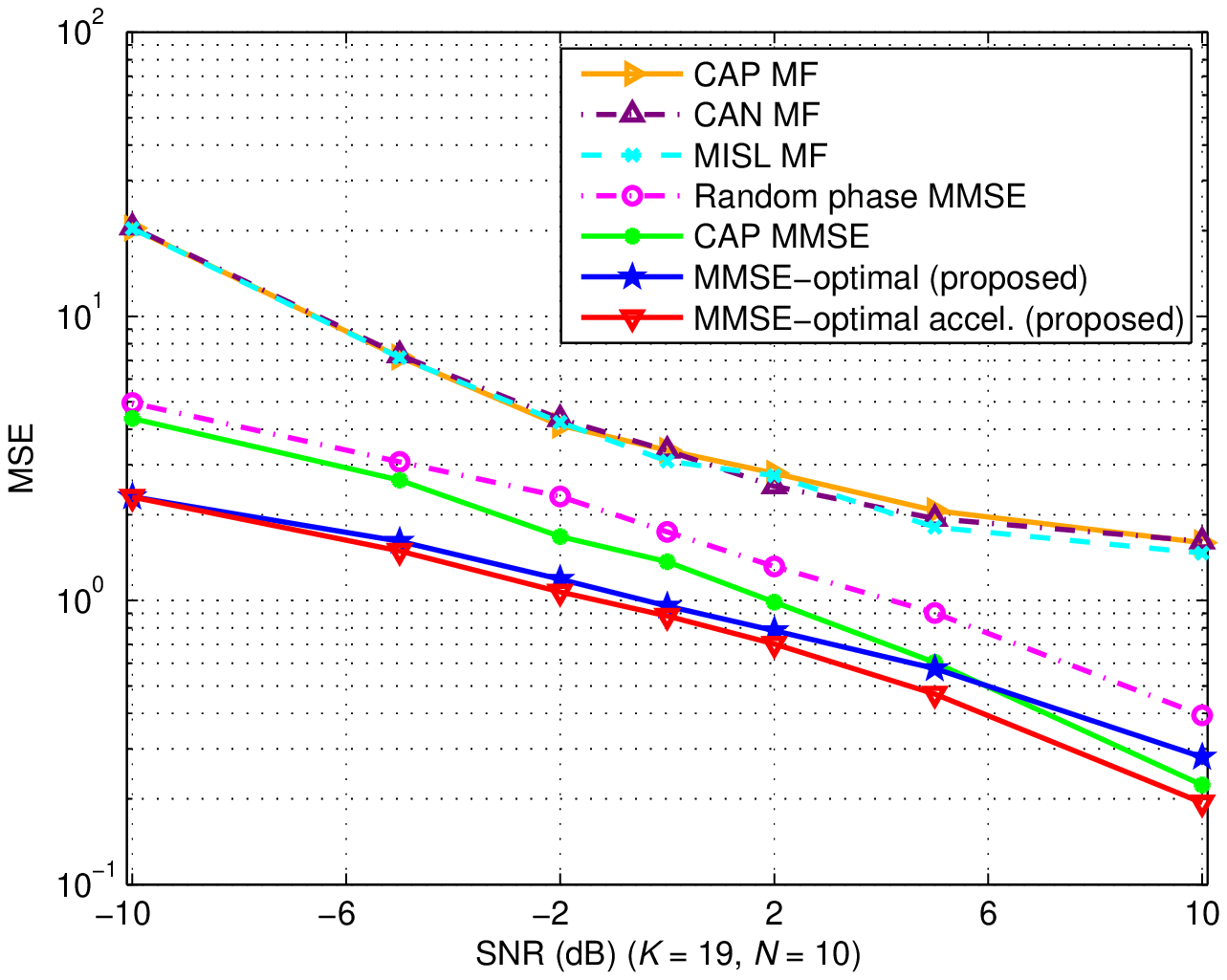}
	\caption{MSE of SISO channel estimates with different unimodular training sequences. The results are averaged over 200 Monte Carlo simulations.}
	\label{figs_unim_MMSE_K19N10}
\end{figure}

\begin{figure}[!t]
	\centering
	\includegraphics[width=3.4in]{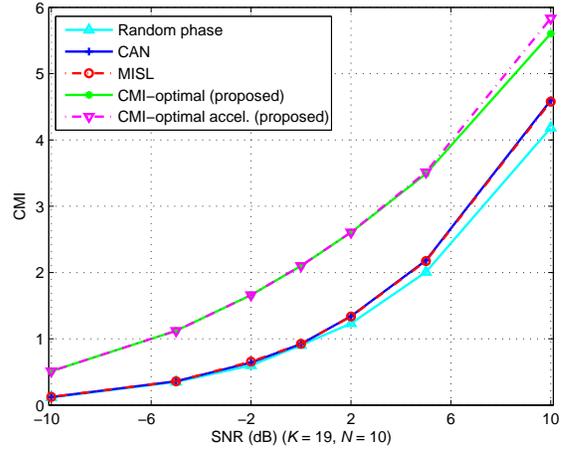}
	\caption{The CMI with different unimodular training sequences for SISO channels. The results are averaged over 200 Monte Carlo simulations.}
	\label{fig_CMI_K29N50}
\end{figure}

\section{Numerical Examples}
\label{sec:simulations}
In this section, we employ proposed algorithms to design unimodular and low PAR sequences for channel estimation. For SISO channels, we compare the channel estimation error and the obtained conditional mutual information of our proposed sequences with that of low sidelobe or random phases. For MIMO channel estimation, the same performance metrics are compared for our proposed sequences, sequences of good auto- and cross-correlation properties, and sequences of random phases. Then we show the advantage of optimal low PAR sequences over sequences of random phases in the MIMO channel estimation.

\subsection{Unimodular Sequences for SISO Channel Estimation}
\label{ssec:SISO_unim_simu}
In this subsection, numerical results are presented to illustrate the advantage of considering the prior information in the design of unimodular sequences for channel estimation and conditional mutual information maximization. Let $N_{t}=N_{r}=1$, and we can apply Algorithm \ref{alg:Optimal_unimodular_sequence} and Algorithm \ref{alg:Accelerated_MM} to design optimal unimodular training sequences for a SISO channel. We compute the MMSE estimates with our proposed sequences, sequences of low sidelobe, and sequences of random phase, and then compare the resulting MSE with matched filtering (MF) using low sidelobe sequences.

The underlying channel impulse response is chosen by $\mathbf{h}_{\mathrm{true}}\sim\mathcal{CN}\left(\mathbf{0}_{K+1},\mathbf{R}_{\mathrm{true}}\right)$
with length $K+1=20$
, and $\left(\mathbf{R}_{\mathrm{true}}\right)_{i,j}=0.9^{\left|i-j\right|}0.9^{\frac{i-1}{2}}0.9^{\frac{j-1}{2}}$ for $i,j=1,\dots,K+1$. The channel is thus correlated with exponentially decreasing power with respect to time delay, which corresponds to the correlated scattering environment with multipath fading in wireless communications \cite{proakisdigital}. The length of training sequence is $N=10$. The channel noise is set to be $\mathbf{v}\sim\mathcal{CN}\left(\mathbf{0}_{N+K},\mathbf{W}\right)$
with $\left(\mathbf{W}\right)_{i,j}=0.2^{\left|i-j\right|}$
for $i,j=1,\dots,N+K$. Considering the inaccuracy of channel covariance matrix in hand, the optimal unimodular sequence $\mathbf{u}$ is designed under the assumed prior $\mathbf{h}_{0}\sim\mathcal{CN}\left(\mathbf{0}_{K+1},\mathbf{R}_{0}\right)$ and $\left(\mathbf{R}_{0}\right)_{i,j}=0.8^{\left|i-j\right|}0.8^{\frac{i-1}{2}}0.8^{\frac{j-1}{2}}$. The mean square error (MSE) of the channel estimator is then
\begin{equation}\label{eq:simu_MMSE}
	\mathrm{MSE}(\hat{\mathbf{h}}_{\mathrm{MMSE}})=\| \hat{\mathbf{h}}_{\mathrm{MMSE}}-\mathbf{h}_{\mathrm{true}}\|^{2}_{2},
\end{equation}
where $\hat{\mathbf{h}}_{\mathrm{MMSE}}$ is given by \eqref{eq:MMSE_estimator} and $\mathbf{S}=\mathcal{T}\left(\mathbf{u}\right)$. Based on the true channel covariance matrix, the conditional mutual information obtained with training sequence $\mathbf{u}$ is
\begin{equation}\label{eq:simu_CMI}
	\mathrm{CMI}\left(\mathbf{u}\right)=\frac{1}{2}\log\det\left(\mathbf{I}+\mathbf{R}_{\mathrm{true}}\mathbf{\mathbf{S}}^{H}\mathbf{W}^{-1}\mathbf{S}\right).
\end{equation}
The signal-to-noise ratio (SNR) is defined as
\begin{equation}\label{eq:simu_SNR}
\mbox{SNR}=10\log_{10}\frac{\left\Vert \mathbf{u}\right\Vert ^{2}\left/N\right.}{\mathrm{Tr}\left(\mathbf{W}\right)\left/\left(N+K\right)\right.}\left(\mathrm{dB}\right).
\end{equation}
For different values of SNR, the resulting MSE and CMI are approximated by running 200 times Monte Carlo simulations. In our simulations, both Algorithm \ref{alg:Optimal_unimodular_sequence} and Algorithm \ref{alg:Accelerated_MM} are initialized with unimodular sequences of random phases uniformly distributed in $[0,2\pi]$. And the algorithms are considered to be converged when the difference between two consecutive updates is no larger than $10^{-6}$, i.e., $\|\mathbf{u}^{(t+1)}-\mathbf{u}^{(t)}\|_{2}\leq 10^{-6}$.

Fig. \ref{figs_unim_MMSE_K19N10} shows the MSE of different channel estimates after training with different unimodular sequences. Both CAP and CAN were proposed to design sequences with low sidelobes, or good correlation properties, and sequences designed by CAP was employed to estimate channel impulse response with the matched filter \cite{stoica2009new}. It was claimed that MISL could further reduce the sidelobes of the designed unimodular sequences \cite{song2014optimal}, with which channel estimate by matched filtering was also compared herein. The resulting MSE of our proposed sequence, MMSE-optimal accel., by the accelerated scheme Algorithm \ref{alg:Accelerated_MM} is lower than that of low sidelobes and that of random phases, especially in the low SNR scenarios. Therefore, the good correlation properties do not guarantee a good channel estimate when the length of the training sequence is limited with respect to the length of the channel impulse response. Note that sequence MMSE-optimal by Algorithm \ref{alg:Optimal_unimodular_sequence} achieves almost the same performance as that of MMSE-optimal accel., but the resulting MSE degrades a little bit in the high SNR case as it needs more iterations to converge. The convergence of Algorithm \ref{alg:Optimal_unimodular_sequence} and Algorithm \ref{alg:Accelerated_MM} will be illustrated in Section \ref{ssub:convergence}.

The obtained CMI for different unimodular sequences are shown in Fig. \ref{fig_CMI_K29N50}. Although by definition \eqref{eq:simu_CMI}, the resulting CMI only depends on the channel statistics without being affected by the channel realizations, Monte Carlo simulations are still conducted for 200 times to avoid the effects from local minima. Expectedly, sequences obtained by CAN and MISL produces almost the same CMI. By incorporating the prior channel information into the sequence design, however, the CMI obtained is improved.

\begin{figure}[!t]
	\centering
	\includegraphics[width=3.4in]{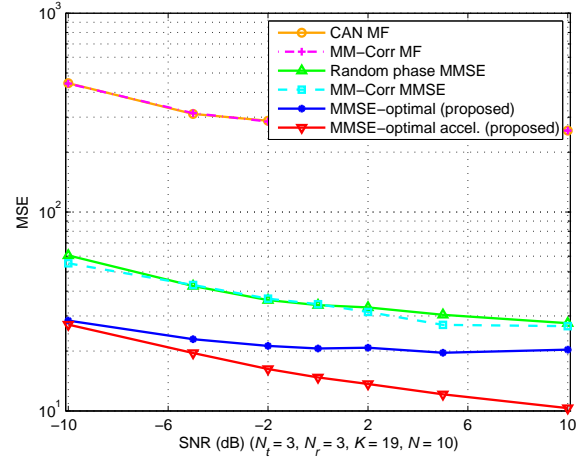}
	\caption{MSE of MIMO channel estimates with different unimodular training sequences. The results are averaged over 100 Monte Carlo simulations.}
	\label{figs_unim_MIMO_MMSE}
\end{figure}

\begin{figure}[!t]
	\centering
	\includegraphics[width=3.4in]{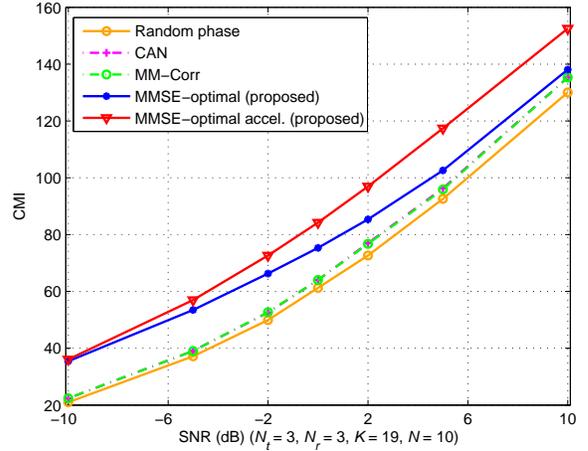}
	\caption{The CMI with different unimodular training sequences for MIMO channels. The results are averaged over 100 Monte Carlo simulations.}
	\label{figs_unim_MIMO_CMI}
\end{figure}

\subsection{Unimodular Sequences for MIMO Channel Estimation}
\label{ssec:MIMO_unim_simu}
In this subsection, we compare the optimal unimodular sequences with those of good correlation properties \cite{song2015seqencce} or random phases for MIMO channels. As in the case of SISO channels, two performance metrics are considered, namely the channel MSE and CMI.

Suppose the MIMO channel has $N_{t}=3$ transmit antennas and $N_{r}=3$ receive antennas, with the length of the channel impulse $K+1=20$. The vectorized channel impulse response $\mathbf{h}_{\mathrm{true}}$ is drawn from a circular complex Gaussian distribution $\mathcal{CN}\left(\mathbf{0}_{N_{t}N_{r}(K+1)},\mathbf{R}_{\mathrm{true}}\right)$. Each channel coefficient $\left(\mathbf{h}_{\mathrm{true}}\right)_{i},i=1,\dots,N_{t}N_{r}(K+1)$ is associated with a triple set $\left(n_{t},n_{r},k\right)$, where $n_{t}=1,\dots,N_{t}$ and $n_{r}=1,\dots,N_{r}$ are indices of transmit and receive antenna, respectively, and $k=0,\dots,K$ is the channel delay. And each entry $\left(\mathbf{R}_{\mathrm{true}}\right)_{i,j}$ of the covariance matrix describes the correlation between the channel coefficient of the triple set $\left(n_{t1},n_{r1},k_{1}\right)$ and $\left(n_{t2},n_{r2},k_{2}\right)$. Without loss of generality, consider
\begin{equation}\label{eq:simu_MIMO_corr}
	\mathbf{R}_{\mathrm{true}} = \mathbf{R}_{r}\otimes \mathbf{R}_{d} \otimes \mathbf{R}_{t}
\end{equation}
where $\left(\mathbf{R}_{r}\right)_{n_{r1},n_{r2}}=\rho_{1}^{\left|n_{r1}-n_{r2}\right|}$ and $\left(\mathbf{R}_{t}\right)_{n_{t1},n_{t2}}=\rho_{3}^{\left|n_{t1}-n_{t2}\right|}$ characterizes, respectively, the correlation between transmit antennas and the correlation between receive antennas, and $\left(\mathbf{R}_{d}\right)_{k_{1},k_{2}}=\rho_{2}^{\left|k_{1}-k_{2}\right|}$ is an exponentially decaying correlation with respect to the channel delay. For the true channel impulse response $\mathbf{h}_{\mathrm{true}}$, we set $\rho_{1}=\rho_{3}=0.9$ and $\rho_{2}=0.7$. In the optimal unimodular training sequence design, the channel prior $\mathbf{h}_{0}$ is assumed to follow a circularly complex Gaussian distribution with zero mean and covariance matrix $\mathbf{R}_{0}$ of the same correlation structure as \eqref{eq:simu_MIMO_corr} and $\rho_{1}=\rho_{3}=0.8$ and $\rho_{2}=0.6$. Each column of noise matrix $\mathbf{V}$ in model \eqref{eq:MIMO_channel_matform} corresponds to a MISO channel, and the vectorized noise is assumed to be colored with a Toeplitz correlation and $\mathrm{vec}\left(\mathbf{V}\right)\sim\mathcal{CN}\left(\mathbf{0}_{(N+K)N_{t}},\mathbf{W}\right)$, with $W_{i,j}=0.2^{|i-j|},i,j=1,\dots,(N+K)N_{r}$. The optimal unimodular training sequences, sequences of good auto- and cross-correlations properties, and sequences of random phases are transmitted and then the corresponding MMSE channel estimators can be obtained. The MSE for each estimate is calculated by \eqref{eq:simu_MMSE} with $\mathbf{S}=\mathcal{T}\left(\mathbf{U}\right)$. The CMI is similarly defined by \eqref{eq:simu_CMI}. The SNR is defined as
\begin{equation}
	\mbox{SNR}=10\log_{10}\frac{\left\Vert \mathbf{U}\right\Vert _{F}^{2}\left/\left(NN_{t}\right)\right.}{\mathrm{Tr}\left(\mathbf{W}\right)\left/\left((N+K)N_{r}\right)\right.}\left(\mathrm{dB}\right).
\end{equation}
The setting for algorithm initialization and convergence are the same as the unimodular case. And the MSE and CMI are averaged over 100 times Monte Carlo simulations for different values of SNR.

Fig. \ref{figs_unim_MIMO_MMSE} shows the MSE of MMSE channel estimates with different unimodular training sequences and SNR's. The length of sequence for each transmit antenna is $N=10$. It is obvious that the optimal unimodular sequences, both MMSE-optimal by Algorithm \ref{alg:Optimal_unimodular_sequence} and MMSE-optimal accel. by Algorithm \ref{alg:Accelerated_MM}, produce smaller MSE than that of random phases or good auto- and cross-correlation properties (Good-Corr). Also notice that there is a gap between two curves of MSE of MMSE-optimal and MMSE-optimal accel. This is because Algorithm \ref{alg:Optimal_unimodular_sequence} needs much more iterations to be converged for MIMO channel training sequence design than that of the SISO case. The convergence properties are shown in Section \ref{ssub:convergence}.

In the CMI maximization for MIMO channels, the performances of different unimodular sequences are shown in Fig. \ref{figs_unim_MIMO_CMI} with $N=10$. For different SNR, the optimal unimodular training sequences can achieve larger CMI than sequences of either random phase or good correlation properties.

\begin{figure}[!t]
	\centering
	\includegraphics[width=3.4in]{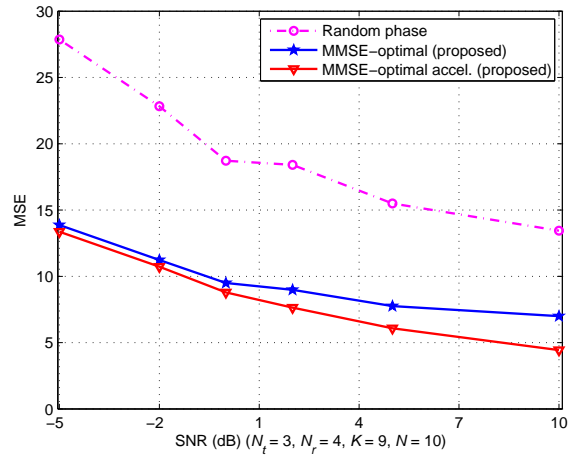}
	\caption{MSE with different low PAR training sequences for MIMO channels. ${\mathrm{PAR}=\{1,2,3\}}$ with power proportions among three antennas: $1:2:3$. The results are averaged over 100 Monte Carlo simulations.}
	\label{figs_lowPAR_MIMO_MMSE}
\end{figure}

\begin{figure}[!t]
	\centering
	\includegraphics[width=3.4in]{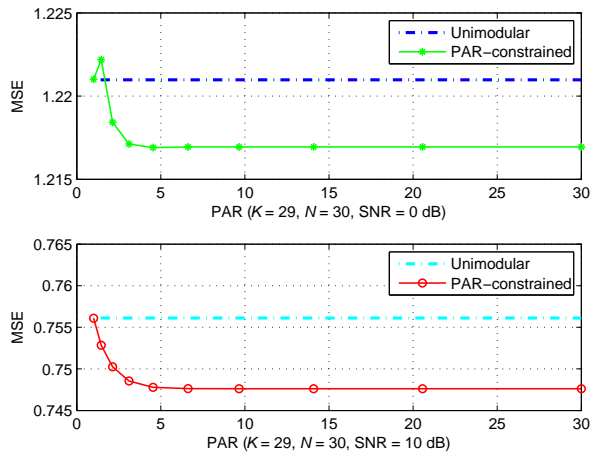}
	\caption{MSE for SISO channel estimation with PAR-constrained sequences or unimodular sequences. The results are averaged over 100 Monte Carlo simulations.}
	\label{figs_SISO_unim_vs_lowPAR_MMSE}
\end{figure}

\subsection{Low PAR Sequences for MIMO Channel Estimation}
Consider the MIMO channel of the same conditions described in Section \ref{ssec:MIMO_unim_simu}. We employ Algorithm \ref{alg:Optimal_PAR_sequence} and its accelerated scheme to design low PAR sequences for the application of MMSE channel estimation. In Fig. \ref{figs_lowPAR_MIMO_MMSE}, MMSE-optimal and MMSE-optimal accel. are obtained by Algorithm \ref{alg:Optimal_PAR_sequence} and its accelerated scheme, respectively. 
It is demonstrated that both optimal training sequences achieve much smaller MSE than low PAR sequences of random phases. Like the results for Algorithm \ref{alg:Optimal_unimodular_sequence} and Algorithm \ref{alg:Accelerated_MM} in the previous subsections, MMSE-optimal renders an larger MSE than MMSE-optimal accel. especially in the high SNR cases. An example of convergence of both algorithms are shown in Section \ref{ssub:convergence}. In Fig. \ref{figs_SISO_unim_vs_lowPAR_MMSE}, we also compare resulting MSE of unimodular sequences and sequences of different values of PAR.

\subsection{Convergence of Proposed Algorithms}
\label{ssub:convergence}
\begin{figure}[!t]
	\centering
	\includegraphics[width=3.4in]{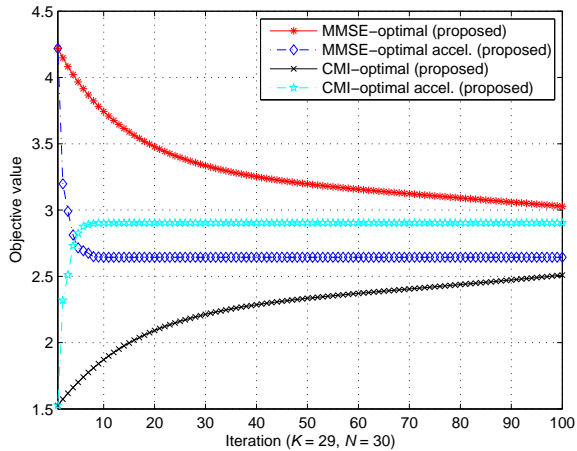}
	\caption{Convergence of algorithms for optimal unimodular sequence design for SISO channel estimation, $\mathrm{SNR}=-5$ dB.}
	\label{fig_SISO_unim_cvg}
\end{figure}

\begin{figure}[!t]
	\centering
	\includegraphics[width=3.4in]{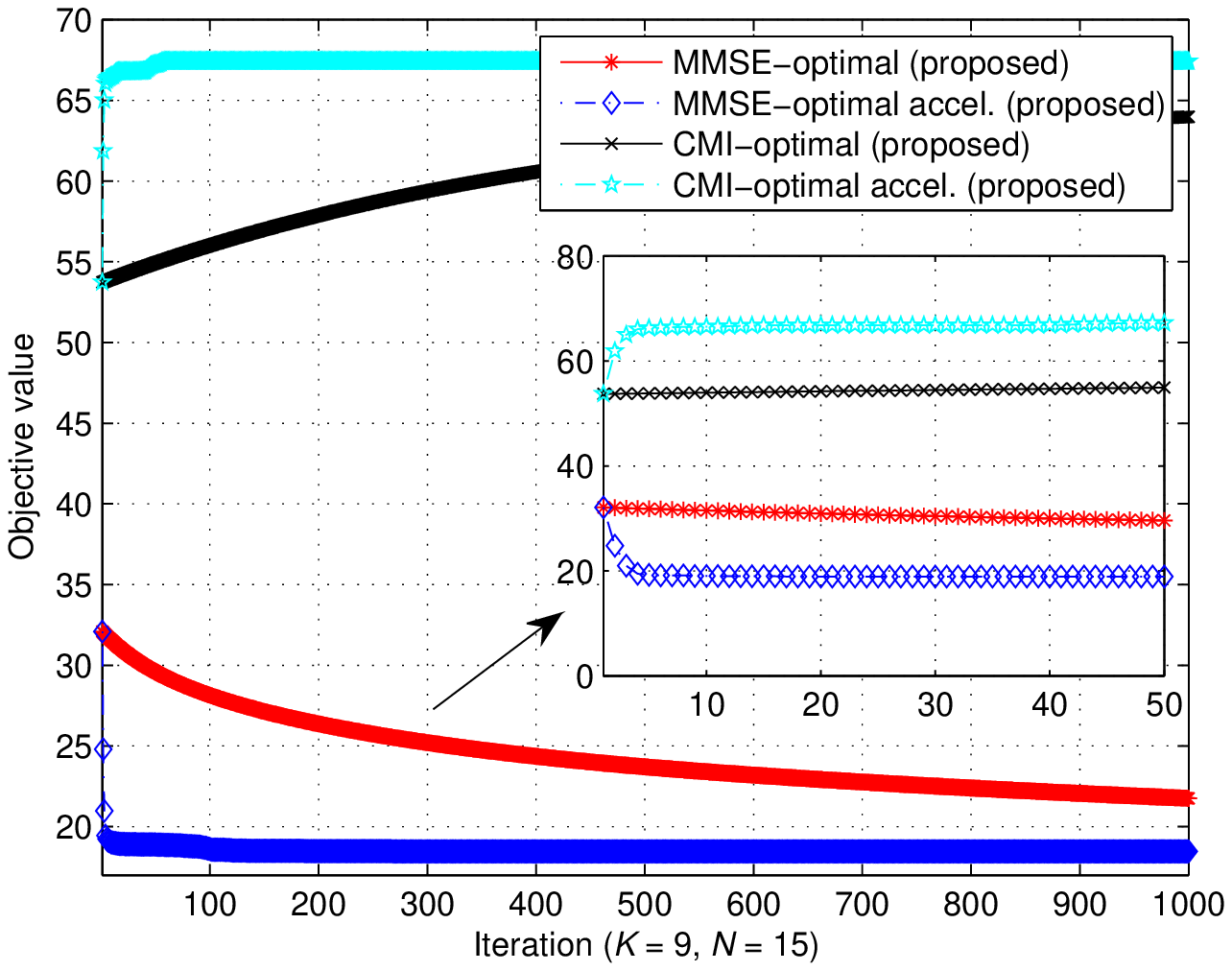}
	\caption{Convergence of algorithms for optimal unimodular sequence design for MIMO channel estimation, $N_{t}=3$, $N_{r}=4$, and $\mathrm{SNR}=-5$ dB.}
	\label{fig_MIMO_unim_cvg}
\end{figure}

\begin{figure}[!t]
	\centering
	\includegraphics[width=3.4in]{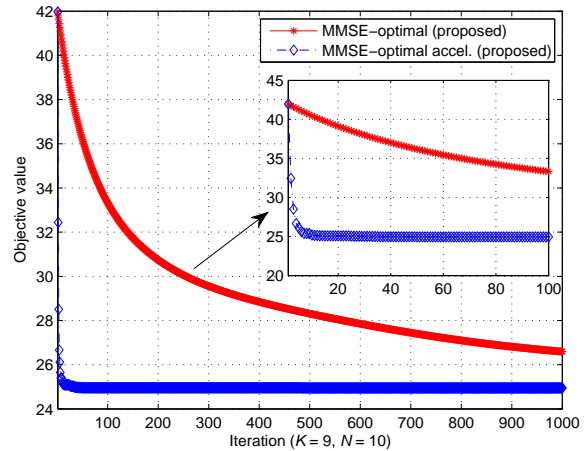}
	\caption{Convergence of algorithms for optimal low PAR sequence design for MIMO channel estimation, $N_{t}=3$, $N_{r}=4$, and $\mathrm{SNR}=-5$ dB.}
	\label{fig_MIMO_lowPAR_cvg}
\end{figure}

Experimental results are given to show the convergence properties of proposed algorithms for the MMSE minimization problem and the CMI maximization problem with unimodular constraints or low PAR constraints. The setting for algorithm initialization and convergence criteria are the same as previous subsections. First, we experiment with Algorithm \ref{alg:Optimal_unimodular_sequence} and Algorithm \ref{alg:Accelerated_MM} for both MMSE minimization and CMI maximization in SISO channel unimodular training sequence design. Fig. \ref{fig_SISO_unim_cvg} shows the objective values with respect to algorithm iterations. In both problems, Algorithm \ref{alg:Optimal_unimodular_sequence} converge monotonically to a stationary point though slowly. With acceleration techniques, however, Algorithm \ref{alg:Accelerated_MM} renders an very fast convergence. The same convergence properties can be seen in Fig. \ref{fig_MIMO_unim_cvg}, where unimodular sequences for MIMO channel estimation are considered with $N_{t}=3$, $N_{r}=4$. Within the same MIMO channel setting, Algorithm \ref{alg:Optimal_PAR_sequence} and its accelerated scheme are applied to design low PAR sequences. The convergence of both algorithms are shown in Fig. \ref{fig_MIMO_lowPAR_cvg}. Note that in those three examples, the algorithms Algorithm \ref{alg:Optimal_unimodular_sequence} and Algorithm \ref{alg:Optimal_PAR_sequence} converge slower than the accelerated scheme especially in designing sequences for MIMO channels with large values of SNR. This is due to successive majorizations or minorizations applied in the derivation of algorithms and thus explains the difference between two training sequences in terms of the resulting MSE and CMI. 

\section{Conclusion}
\label{sec:conclusion}
In this paper, optimal training sequences with unimodular constraint and low PAR constraints are considered. The optimal sequence design problem is formulated by minimizing the MMSE criterion and maximizing the CMI criterion. The formulated problems are nonconvex and efficient algorithms are developed based on the majorization-minimization framework. Furthermore, the acceleration scheme is derived using the SQUAREM method. All the proposed algorithms are guaranteed to monotonically converge to a stationary point. Numerical results show that the optimal unimodular sequences can improve either the accuracy of channel estimate or the CMI compared with those of sequences with good correlation properties or random phases. Under the same criteria, the optimal sequence design with low PAR constraint is also studied, for which the similar algorithms to unimodular case are derived. Numerical examples show that the optimal low PAR sequences perform better than that of random phases.




\ifCLASSOPTIONcaptionsoff
  \newpage
\fi


\bibliographystyle{IEEEtran}

%


%
%




\end{document}